\documentclass[final,5p,times,twocolumn]{elsarticle}

\usepackage{amsmath}
\usepackage{amssymb,bm} 
\usepackage{amsthm,mathtools}
\usepackage{theoremref}
\usepackage[shortlabels]{enumitem}
\usepackage{color}
\usepackage{apptools}
\usepackage{graphicx}
\usepackage{soul}
\usepackage{subcaption}
\newtheorem{thm}{Theorem}
\newtheorem{lem}{Lemma}
\newtheorem{prop}{Proposition}
\newtheorem{cor}{Corollary}
\newtheorem{defn}{Definition}
\newtheorem{asmp}{Assumption}

\newtheorem{exmp}{Example}
\newtheorem{claim}{Claim}
\newtheorem*{exmp*}{Example}
\newtheorem{rem}{Remark}

\newcommand*{\QEDB}{\hfill\ensuremath{\square}}%
\newcommand{\D}[1]{\text{D}_{#1}}

\begin{document}

\begin{frontmatter}



\vspace{-1cm}
\title{\vspace{-0.5cm}{Initialization-Free Lie-Bracket Extremum Seeking in $\mathbb{R}^n$}\tnoteref{titlelabel}}

\tnotetext[titlelabel]{This work was supported in part by the grants NSF ECCS CAREER 2305756 and AFOSR YIP: FA9550-22-1-0211.}

\author[a]{Mahmoud Abdelgalil}
\ead{mabdelgalil@ucsd.edu}
\author[a]{Jorge Poveda}
\ead{jipoveda@ucsd.edu}
\affiliation[a]{organization={Department of Electrical and Computer Engineering, University of California, San Diego},
            city={La Jolla, 92093},
            state={California},
            country={USA}}




\begin{abstract}
Stability results for extremum seeking control in $\mathbb{R}^n$ have predominantly been restricted to local or, at best, semi-global practical stability. Extending semi-global stability results of extremum-seeking systems to unbounded sets of initial conditions often demands a stringent global Lipschitz condition on the cost function, which is rarely satisfied by practical applications. In this paper, we address this challenge by leveraging tools from higher-order averaging theory. In particular, we establish a novel second-order averaging result with \emph{global} (practical) stability implications. By leveraging this result, we characterize sufficient conditions on cost functions under which uniform global practical asymptotic stability can be established for a class of extremum-seeking systems acting on static maps. Our sufficient conditions include the case when the gradient of the cost function, rather than the cost function itself, satisfies a global Lipschitz condition, which covers quadratic cost functions. Our results are also applicable to vector fields that are not necessarily Lipschitz continuous at the origin, opening the door to non-smooth Lie-bracket ES dynamics. We illustrate all our results via different analytical and/or numerical examples.
\end{abstract}


\end{frontmatter}


\section{Introduction}
%
%
Extremum seeking (ES) systems are some of the most popular real-time derivative-free optimization algorithms developed during the last century \cite{Leblanc1922}. The stability and robustness guarantees, simplicity of implementation, and model-agnostic nature of ES make it an attractive option for a large number of practical model-free control problems \cite{ariyur2003real,tan2010extremum,scheinker2024100}. The classical tool for analyzing the stability of ES systems is a (first-order) combination of singular perturbation and averaging theory that enables local practical stability results \cite{krstic2000stability}. These ideas have also been used to study ES systems that approximate Newton's method rather than gradient descent \cite{nesic2010unifying,ghaffari2012multivariable}. In addition, it has been shown that standard ES systems also have non-local stability properties \cite{tan2006non,tan2009global,mimmo2022uniform}, and different schemes have been developed to tackle control and optimization problems in more general systems, such as systems with delays \cite{oliveira2016extremum}, partial differential equations \cite{oliveira2020multivariable}, as well as hybrid dynamical systems \cite{poveda2017framework,poveda2021robust,abdelgalil2023multi}.


On the other hand, recent years have seen the emergence of an alternative approach for designing and analyzing the stability of ES systems based on higher-order averaging \cite{durr2013lie,abdelgalil2022recursive}. Employing higher-order averaging techniques can offer some flexibility in the design and analysis of the exploration-exploitation mechanism of the system, which led to the discovery of ES algorithms with desirable properties such as bounded update rates \cite{scheinker2014extremum}, vanishing amplitudes \cite{scheinker2014non,abdelgalil2021lie}, and even local exponential/asymptotic stability \cite{grushkovskaya2018class}. Higher-order averaging is also better suited for geometric settings when ES is performed on manifolds \cite{durr2013examples,durr2014extremum,abdelgalil2023multi,abdelgalil2023lie}, or when additional structure is imposed on the exploration dynamics, e.g. when the exploration is done through a Levi-Civita connection associated with a mechanical system \cite{suttner2022extremum,suttner2023extremum,suttner2023attitude,suttner2023nonholonomic}.

Irrespective of the nature of the averaging tool used for the analysis and design of the algorithm, when the extremum seeking problem is defined on smooth compact boundaryless manifolds, achieving uniform global stability results (either practical or asymptotic) is in general not possible due to topological obstructions that apply to continuous-time systems evolving on such sets \cite{jongeneel2023topological}, a limitation that also applies to time-varying periodic systems, see \cite[Sec. 4.1]{mayhew2011topological}. However, when the ES problem is defined in $\mathbb{R}^n$, such obstructions do not emerge, and, in principle, it might be possible to achieve global extremum seeking. Nevertheless, the majority of results on ES in $\mathbb{R}^n$ have achieved, at best, semi-global practical asymptotic stability \cite{nesic2010unifying,tan2006non,durr2013lie,abdelgalil2022recursive,scheinker2014extremum}. Such results enable convergence from arbitrarily large pre-defined compact sets of initial conditions by appropriately tuning the parameters of the controller. However, without further re-tuning of these parameters, solutions initialized (or pushed via perturbations) outside of these pre-defined compact sets might exhibit finite escape times. Recently, global practical convergence properties were studied in \cite{mimmo2021extremum} using a normalized scheme, and also in \cite{lauand2023quasi} using tools from quasi-stochastic approximation theory. Nevertheless, results that assert \emph{uniform global practical asymptotic stability} (characterized by, e.g., $\mathcal{KL}$ bounds) in ES controllers are still absent in the literature. Indeed, one of the main limitations that makes achieving such a result difficult using standard averaging theory stems from the global Lipschitz conditions that are usually required in the vector fields of the dynamics, see \cite[Ch.10]{khalil2002nonlinear}, \cite[Sec. 6.1]{abdelgalil2023multi}, a condition that is violated even in the simplest ES problems corresponding to cost functions characterized by quadratic maps.
%

%
Based on the above background, the main contribution of this paper is to establish that certain ES systems can achieve uniform global (practical) stability results once one moves from \emph{first-order} averaging-based feedback designs, such as those considered in \cite{mimmo2021extremum,tan2006non,krstic2000stability,poveda2017framework}, to \emph{second-order} averaging-based feedback designs, similar in spirit to those studied in \cite{durr2013examples,durr2014extremum,abdelgalil2023multi,abdelgalil2023lie,suttner2022extremum,suttner2023extremum}, but using a different averaging tool for the purpose of analysis. In particular, the main contribution of this paper is twofold: First, we introduce a novel \emph{second-order} averaging theorem with \emph{global} practical stability implications for a general class of highly-oscillatory systems under appropriate assumptions on the maps involved. For standard (i.e., first-order) averaging, global stability results have been studied in \cite{sastry1990adaptive} for ODEs, and in \cite{abdelgalil2023multi} for hybrid systems. However, to the best of our knowledge, a result of this nature was absent in the literature of second-order averaging. Furthermore, in contrast to existing results on second-order averaging \cite{durr2013lie}, the novel averaging tool enables the relaxation of the local Lipschitz condition on the vector field at the origin to mere continuity. This relaxation opens the door to novel non-smooth dynamics that could potentially lead to improved transient performance away from the origin. Second, we use the aforementioned averaging results to establish uniform global practical asymptotic stability properties for a class of ES systems for which a variety of ``typical'' cost functions apply, including quadratic maps, and, more generally, strongly convex functions with smooth gradients. However, we also show that convexity of the cost function is, in general, not a necessary condition to achieve global ES under the algorithms studied in this paper. Different analytical and numerical examples are presented to illustrate our results.

\vspace{0.1cm}
%
The rest of the manuscript is organized as follows. We begin by introducing our notation in Section \ref{sec:notation}. Our averaging results are then presented in Section \ref{sec:averaging_results}. In Section \ref{sec:es_section}, we apply the results of Section \ref{sec:averaging_results} to study a class of extremum seeking systems that attain global (practical) stability properties. The proofs of the results are presented in Section \ref{sec:proofs}. Finally, the conclusions and future work are discussed in Section \ref{sec:conclusions}.
\section{ Preliminaries}\label{sec:notation}
\subsection{Notation}

 We use $\mathbb{R}_{\geq0}$ to denote the set of non-negative real numbers and $\mathbb{R}_{>0}$ to denote the set of positive real numbers. Similarly, we use $\mathbb{Q}_{>0}$ to denote the set of positive rational numbers and $\mathbb{N}_{\geq 1}$ to denote the set of positive integers. The 2-norm of a vector $x\in\mathbb{R}^n$ is denoted by $|x|:=\sqrt{x^\top x}$, and the operator $2$-norm of a matrix $A\in\mathbb{R}^{m\times n}$ is also denoted as $|A|:=\sup\{|Ax|: x\in\mathbb{R}^n,\,|x|=1\}$. If $x\in\mathbb{R}^n$ and $y\in\mathbb{R}^m$ are vectors, we use $(x,y)=[x^\top,y^\top]^\top\in\mathbb{R}^{n+m}$ to denote the vector consisting of their concatenation. Given functions $f:\mathbb{R}^n\rightarrow \mathbb{R}^m$ and $g:\mathbb{R}^m\rightarrow \mathbb{R}^l$, we use $g\circ f:\mathbb{R}^n\rightarrow\mathbb{R}^l$ to denote their composition, i.e. $g\circ f(x) = g(f(x))$. We use $\mathcal{C}^0$ to denote the class of continuous functions, and $\mathcal{C}^k$ to denote the class of functions that are $k$-times continuously differentiable, for $k\geq1$. Given a closed set $K\subset\mathbb{R}^n$, the function $f$ is said to be $\mathcal{C}^k$ on $K$ if there exists an open neighborhood $\mathcal{U}\subset\mathbb{R}^n$ such that $K\subset\mathcal{U}$ and $f$ is $\mathcal{C}^k$ on $\mathcal{U}$. For each $\delta\in\mathbb{R}_{>0}$, we denote the closed ball of radius $\delta$, centered at the origin, by $\delta\mathbb{B}$, i.e. $\delta\mathbb{B}:=\{x\in\mathbb{R}^n:\,|x| \leq \delta\}$. Given a set $\mathcal{A}\subset\mathbb{R}^n$, we use $\text{cl}(\mathcal{A})$ to denote the closure of $\mathcal{A}$ with respect to the natural topology in $\mathbb{R}^n$. When $f\in\mathcal{C}^1$ is a vector-valued map, $\D{} f$ denotes the Jacobian of $f$. If $f\in\mathcal{C}^2$ is a real-valued function, then $\nabla f$ denotes the gradient of $f$, and $\nabla^2 f$ is the Hessian of $f$, i.e. $\nabla^2 f = \D{}(\nabla f)$. If $f\in\mathcal{C}^1$ and $f=f(x_1,\dots,x_n)$ is vector-valued, then $\D{x_i} f$ denotes the Jacobian of $f$ with respect to the $i^{th}$ argument. The map $\pi_{i}:\mathbb{R}^{n_1}\times\cdots\times\mathbb{R}^{n_k}\rightarrow \mathbb{R}^{n_i}$ is the canonical projection onto the $x_i$-factor, which is defined by $\pi_{i}(x_1,\dots,x_k)=x_i$. A class $\mathcal{K}$-function is a strictly increasing continuous function $\alpha:\mathbb{R}_{\geq0}\to\mathbb{R}_{\geq0}$ such that $\alpha(0)=0$. A class $\mathcal{K}_\infty$-function is a class $\mathcal{K}$-function with the additional requirement that $\lim_{\rho\rightarrow +\infty}\alpha(\rho)=+\infty$. A class $\mathcal{KL}$-function $\beta:\mathbb{R}_{\geq0}\times\mathbb{R}_{\geq0}\to\mathbb{R}_{\geq0}$ is a continuous function such that, for every $s\in\mathbb{R}_{\geq0}$, the function $\beta(\cdot,s)$ is a class $\mathcal{K}_\infty$-function, and, for every $r\in\mathbb{R}_{\geq0}$, the function $\beta(r,\cdot)$ is a strictly decreasing function and $\lim_{s\rightarrow +\infty}\beta(r,s) = 0$. To simplify notation, given two (or more) vectors $x_1\in\mathbb{R}^{n_1},x_2\in\mathbb{R}^{n_2}$, we use $(x_1,x_2)\in\mathbb{R}^{n_1+n_2}$ to denote the concatenation of $x_1$ and $x_2$. 
\subsection{Dynamical Systems and Stability Notions}
In this paper, we study continuous-time dynamical systems with states $(x,\tau)\in\mathbb{R}^n\times\mathbb{R}_{\geq 0}$, and dynamics
\begin{align}\label{ODE01}
        \dot{x}= f_\varepsilon(x,\tau),~~~~~~\dot{\tau}&= \varepsilon^{-2},
\end{align}
where $f_{\varepsilon}:\mathbb{R}^n\times\mathbb{R}_{\geq0}\to\mathbb{R}^n$ is a continuous function parameterized by a small constant $\varepsilon>0$. Systems of the form \eqref{ODE01} can model highly oscillatory systems that showcase fast variations of $\tau$ compared to the state $x$. For completeness, the notion of solutions to systems of the form \eqref{ODE01} is reviewed below.
\begin{defn}\label{defn:solution_concept}\normalfont
    For $(x_0,\tau_0)\in\mathbb{R}^n\times\mathbb{R}_{\geq 0}$, a function $(x,\tau):\text{dom}(x,\tau)\rightarrow \mathbb{R}^{n}\times\mathbb{R}_{\geq0}$ is said to be a solution to \eqref{ODE01} from the initial condition $(x_0,\tau_0)$ if: i) there exist $t_s\in\mathbb{R}_{>0}\cup\{\infty\}$ such that $\text{dom}(x,\tau)=[0,t_s)$, ii) $(x(0),\tau(0))=(x_0,\tau_0)$, and iii) the function $(x,\tau)$ is continuously differentiable on $\text{dom}(x,\tau)$ and satisfies
    \begin{align*}
        \frac{dx(t)}{dt}&= f_\varepsilon(x(t),\tau(t)), & \frac{d\tau(t)}{dt}&= \varepsilon^{-2}.
    \end{align*}
    for all $t\in \text{dom}(x,\tau)$.
    A solution $(x,\tau)$ to system \eqref{ODE01} is said to be complete if $t_s=\infty$.  \QEDB 
\end{defn}
Under the action of a $\mathcal{C}^1$ diffeomorphism $\Psi:\mathbb{R}^n\times\mathbb{R}_{\geq 0}\rightarrow \mathbb{R}^n\times\mathbb{R}_{\geq 0}$, a solution $(x,\tau)$ of the ODE \eqref{ODE01} is transformed into a new function $\Psi\circ(x,\tau):\text{dom}(x,\tau)\rightarrow\mathbb{R}^n\times\mathbb{R}_{\geq0}$. The following definition (see, e.g., \cite[p. 183]{lee2012smooth}) characterizes an ODE for which $\Psi\circ(x,\tau)$ is a solution.
\begin{defn}\label{defn:diffeomorphism_action_on_ODE}\normalfont
    Let $\Psi:\mathbb{R}^n\times\mathbb{R}_{\geq 0}\rightarrow \mathbb{R}^n\times\mathbb{R}_{\geq 0}$ be a $\mathcal{C}^1$ diffeomorphism such that $\pi_2\circ\Psi(x,\tau)=\tau$. The \emph{pushforward} of the ODE \eqref{ODE01} under the action of $\Psi$ is the ODE
     \begin{subequations}\label{pushforward}
    \begin{align}\label{eq:perturbed_avged_system_timevaring}
        \dot{x}&= \Psi_*f_\varepsilon(x,\tau), & \dot{\tau}&=\varepsilon^{-2},
    \end{align}
    where the map $\Psi_*f_\varepsilon$ is given by
    \begin{align}\label{eq:push_forward_map_definition}    
    \Psi_*f_\varepsilon&=\left(\D{x}\left(\pi_1\circ\Psi\right)\circ\Psi^{-1}\right)\,f_\varepsilon\circ\Psi^{-1} + \D{\tau}\left(\pi_1\circ\Psi\right)\circ\Psi^{-1}\varepsilon^{-2},
    \end{align}
    defined for all $(x,\tau)\in\mathbb{R}^n\times\mathbb{R}_{\geq0}$. \QEDB 
     \end{subequations}
\end{defn}
To study the (uniform) stability properties of the parameter-dependent system \eqref{ODE01}, we will use the following standard notions (see, e.g., \cite{teel2010averaging}). Note that we do not insist on uniqueness of solutions, but rather impose the appropriate bound (and the property of completeness) to every solution of the system.
\begin{defn}\normalfont
    The origin $x^*=0$ is said to be \emph{uniformly globally practically asymptotically stable (\textbf{UGpAS})} for system \eqref{ODE01} if there exists a class $\mathcal{KL}$-function $\beta$ such that, for every $\nu\in\mathbb{R}_{>0}$, there exists $\varepsilon^*>0$, such that, for all $\varepsilon\in(0,\varepsilon^*)$, each solution $(x,\tau)$ to system \eqref{ODE01} from the initial condition $(x_0,\tau_0)\in \mathbb{R}^n\times\mathbb{R}_{\geq 0}$ satisfies
    \begin{align}\label{KLbound1001}
        |x(t)|&\leq \beta(|x_0|,t) + \nu,
    \end{align}
    for all $t\geq0$. When $v=0$, the origin $x^*=0$ is said to be \emph{uniformly globally asymptotically stable (\textbf{UGAS})} for  \eqref{ODE01}. \QEDB
\end{defn}
When the residual upper-bound $\nu$ in \eqref{KLbound1001} cannot be controlled by the parameter $\varepsilon$, we will study the following property.
\begin{defn}\normalfont
    The origin $x^*=0$ is said to be $\Delta$-\emph{uniformly globally ultimately bounded} ($\Delta$-\textbf{UGUB}) for system \eqref{ODE01} if there exists $\Delta>0$, $\beta\in\mathcal{KL}$, and $\varepsilon^*\in\mathbb{R}_{>0}$, such that for all $\varepsilon\in(0,\varepsilon^*)$, each solution $(x,\tau)$ to system \eqref{ODE01} from the initial condition $(x_0,\tau_0)\in \mathbb{R}^n\times\mathbb{R}_{\geq 0}$ satisfies
    \begin{align}\label{UUbound}
        |x(t)|&\leq \beta(|x_0|,t) + \Delta, 
    \end{align}
    for all $t\geq0$. \QEDB
\end{defn}

\section{On Global Stability via Second-Order Averaging}\label{sec:averaging_results}
\begin{figure}[t]
    \centering
    \includegraphics[width=0.95\linewidth]{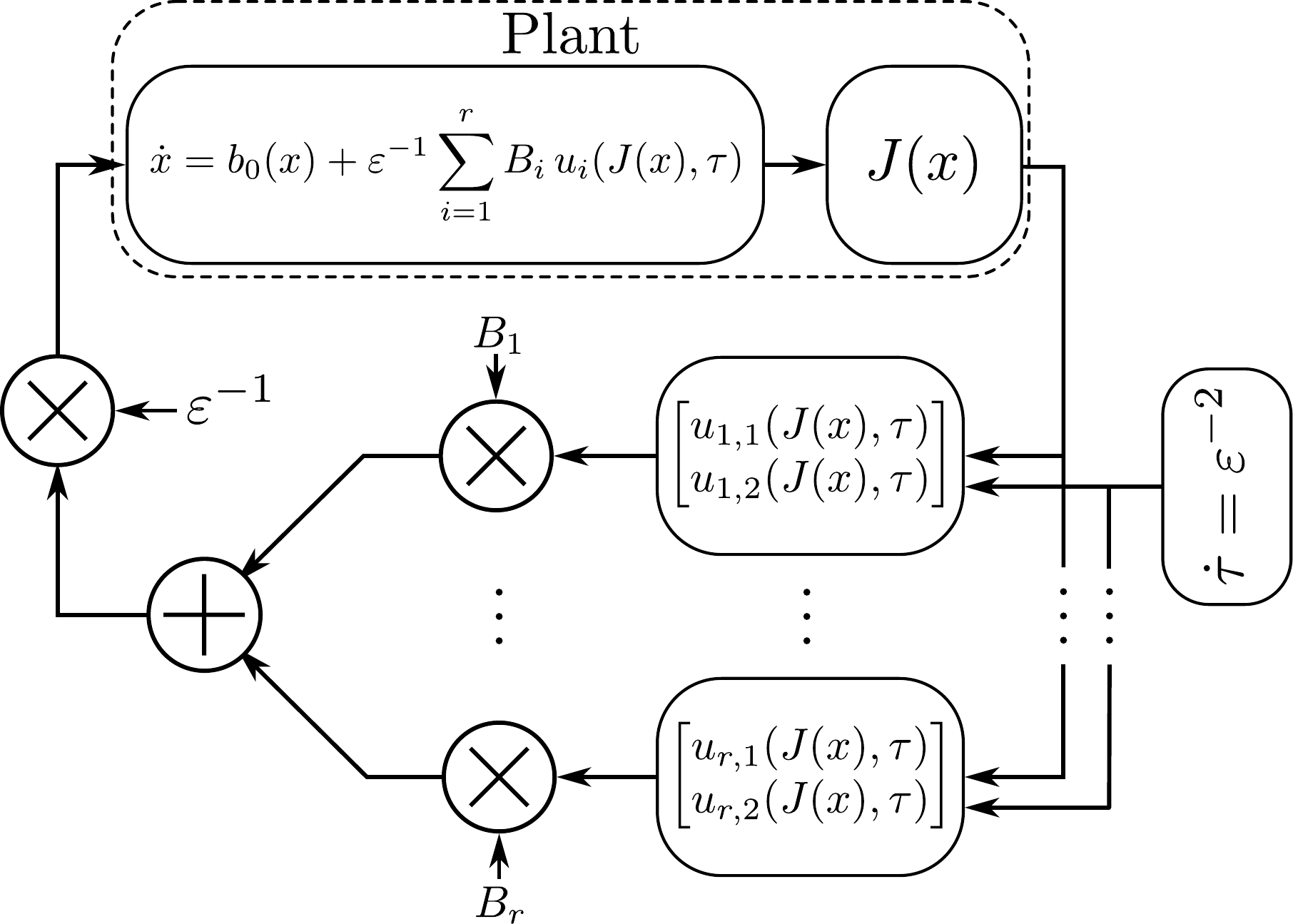}
    \caption[]{Block diagram description of system \eqref{eq:es_system_open_loop}. In the diagram, the matrix $B_i=[b_{i,1},\,b_{i,2}]$ multiplies the vector $u_i(J(x),\tau)=(u_{i,1}(J(x),\tau),\,u_{i,2}(J(x),\tau))$.}
    \label{fig:es_block_diagram}
\end{figure}

We consider a sub-class of systems of the form \eqref{ODE01}, given by
%
\begin{equation}\label{eq:orig_system_timevaring}
        \dot{x}= f_\varepsilon(x,\tau)=\varepsilon^{-1}f_1(x,\tau) + f_2(x,\tau),~~~~~~~\dot{\tau}=\varepsilon^{-2},
\end{equation}
%
where $f_k:\mathbb{R}^n\times\mathbb{R}_{\geq0}\to\mathbb{R}^n$, $k\in\{1,2\}$, are continuous functions, and $\varepsilon>0$. Such types of systems commonly emerge in extremum-seeking \cite{abdelgalil2021lie,durr2013lie} and vibrational control \cite{scheinker2013model}, and they are typically studied via averaging theory. A representative example is given by control-affine systems of the form
\begin{align}\label{eq:es_system_open_loop}
        \dot{x}&=\varepsilon^{-1}\left(\sum_{\substack{i=1}}^{r}\sum_{j=1}^2 b_{i,j}u_{i,j}\left(J(x),\tau\right)\right)+b_0(x), & \dot{\tau}&= \varepsilon^{-2},
\end{align}
where $x\in\mathbb{R}^{n}$, $r\in \mathbb{N}_{\geq \frac{n}{2}}$, $J$ is an application-dependent $\mathcal{C}^2$ cost function to be minimized, $b_{i,j}$ are suitable vectors, $u_{i,j}(\cdot,\cdot)$ is a scalar-valued feedback law to be designed, and $\varepsilon>0$ is a small tunable parameter, see Figure \ref{fig:es_block_diagram} for a block representation of these systems. Particular examples of functions $b_0,b_{i,j},u_{i,j}$ and $J$ will be discussed later in Section \ref{sec:applicationsES}. 

%

%
\subsection{A Global Practical Near-Identity Transformation}
Traditionally, the averaging-based analysis of oscillatory systems relies on the construction of a suitable (first-order) ``near-identity'' transformation that maps the original dynamics into a perturbed version of the so-called average dynamics, see \cite[Ch.10]{khalil2002nonlinear}. Therefore, to study the global stability properties of \eqref{eq:orig_system_timevaring}, we first construct a similar ``second-order'' near-identity transformation, of global nature, and we show how to use this transformation to transform \eqref{eq:orig_system_timevaring} into a perturbed version of its average dynamics.

\vspace{0.1cm}
We consider the following regularity conditions on $f_k$:
\begin{asmp}\thlabel{asmp:vector_fields_condition_time_varying}\normalfont
    There exists $\delta_1\in[0,\infty)$ such that, for all $k\in\{1,2\}$, the following conditions hold
\begin{enumerate}[(a),topsep=0pt,itemsep=-1ex,partopsep=1ex,parsep=1ex]
    \item\label{asmp:vector_fields_condition_time_varying_1} The map $f_k$ is $\mathcal{C}^0$ in $\mathbb{R}^n\times\mathbb{R}_{\geq0}$, and there exist positive constants $L_{k} $ such that
    \begin{gather*}
        |f_k(x_1,\tau)-f_k(x_2,\tau)|\leq L_k |x_1-x_2|,
    \end{gather*}
    for all $x_1,x_2\in\{x\in\mathbb{R}^n: |x|\geq \delta_1\}$ and all $\tau\in\mathbb{R}_{\geq 0}$.
    \item \label{asmp:vector_fields_condition_time_varying_2}There exists $T\in\mathbb{R}_{>0}$ such that
    \begin{align*}
        f_k(x,\tau+T)&=f_k(x,\tau), &  \int_{0}^T f_1(x,\tau) d\tau = 0, 
    \end{align*}
    for all $(x,\tau)\in\mathbb{R}^n\times\mathbb{R}_{\geq 0}$.
    \item\label{asmp:vector_fields_condition_time_varying_3} The map $f_k$ is $\mathcal{C}^{3-k}$ with respect to $x$ in the domain $\{x\in\mathbb{R}^n: |x|\geq \delta_1\}$.
    \item\label{asmp:vector_fields_condition_time_varying_4} There exists $L_{3}>0$ such that
    \begin{align*}
        \hspace{-0.2cm}|\D{x} f_1(x_1,\tau_1) f_1(x_1,\tau_2)-\D{x} &f_1(x_2,\tau_1) f_1(x_2,\tau_2)|\leq L_{3} |x_1-x_2|.
    \end{align*}
    for all $x_1,x_2\in\{x\in\mathbb{R}^n: |x|\geq \delta_1\}$ and all $\tau_1,\tau_2\in\mathbb{R}_{\geq 0}$.
\end{enumerate}
\vspace{-0.3cm}\QEDB
\end{asmp}
In \thref{asmp:vector_fields_condition_time_varying}, the case $\delta_1=0$ is included. However, we allow for positive values of $\delta_1$ to account for the situation in which the regularity of the maps $f_k$ drops from being $\mathcal{C}^{3-k}$ to merely $\mathcal{C}^0$, as required by item \ref{asmp:vector_fields_condition_time_varying_1} in \thref{asmp:vector_fields_condition_time_varying}, near the origin. Relaxing the smoothness requirement allows us to include in our analysis certain non-smooth ES dynamics that have been show to exhibit suitable local exponential/asymptotic stability \cite{scheinker2014non,suttner2017exponential}. 

\vspace{0.1cm}
Next, we introduce the auxiliary functions $\chi_j:\mathbb{R}\rightarrow \mathbb{R}_{\geq 0}$, given by
\begin{align*}
    \chi_1(r)&:= \begin{cases}
        \exp\left(- r^{-1}\right) & r>0\\
        0 & r\leq 0
    \end{cases}, & ~~\chi_2(r)&:= \frac{\chi_1(r)}{\chi_1(r) + \chi_1(1-r)}.
\end{align*}
Also, for each $\epsilon>0$ we consider a vector of constants given by $\delta:=(\delta_1,\delta_2,\delta_3)$, which satisfy
\begin{equation}\label{eq:deltas}
\delta_1\geq0,~~~~~\delta_2\geq(1+\epsilon)\delta_1,~~~~~\delta_3>(1+\epsilon)\delta_2.
\end{equation}
Using the function $\chi_2$ and the vector $\delta$, we define the smooth ``reverse" bump function $\varphi:\mathbb{R}^n\rightarrow [0,1]$ as:
\begin{subequations}\label{eq:smooth_vector_fields}
    \begin{align}\label{bumpfunction}
    \varphi(x)&= \begin{cases}
        \chi_2\left(\frac{|x|-\delta_1}{\delta_2-\delta_1}\right) & \delta_2>\delta_1\\
        1 & \delta_2=\delta_1=0.
    \end{cases}
\end{align}
The function $\varphi$ will be used only for the purpose of analysis, and any similarly defined smooth ``reverse'' bump function suffices for our purposes. The following Lemma states some useful properties of  $\varphi$.
\begin{lem}\thlabel{lem:hole_function_properties}
    Let $\delta_2> \delta_1$. Then, the function $\varphi$ is $\mathcal{C}^\infty$ on $\mathbb{R}^n$, all of its derivatives have the compact support $[\delta_1,\delta_2]$, and it satisfies:
    \begin{enumerate}[(a),topsep=0pt,itemsep=-1ex,partopsep=1ex,parsep=1ex]
    \item $\varphi(x)=1$ for all $x\in\{x'\in\mathbb{R}^n:\, |x|\geq \delta_2\}$.
    \item $\varphi(x)=0$ for all $x\in\{x'\in\mathbb{R}^n:\, |x|\leq \delta_1\}$.
    \end{enumerate}
\end{lem}
\noindent\textbf{Proof: } Follows by \cite[Lemmas 2.20-2.22]{lee2012smooth} and the construction of the argument of $\chi_2$. \QEDB 

\begin{figure}[t!]
    \centering
    \includegraphics[width=0.7\linewidth]{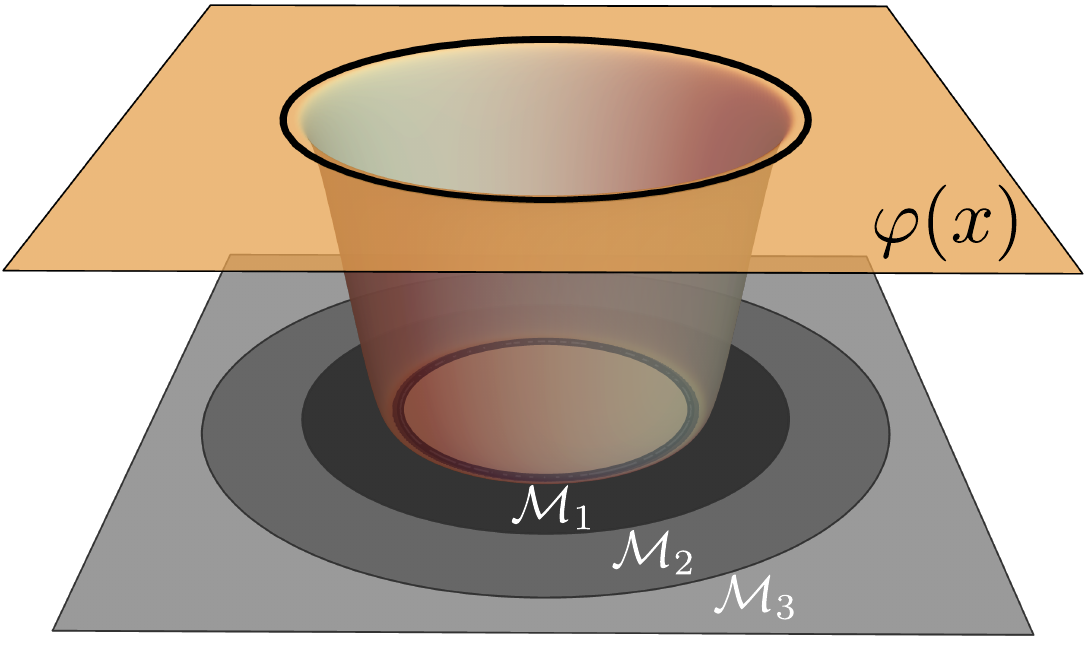}
    \caption[]{Visual depiction of $\varphi$ and the sets $\mathcal{M}_j$ for $j\in\{1,2,3\}$.}
    \label{fig:bump_function}
\end{figure}
\vspace{0.1cm}
To state our first result, and using $\varphi$, we introduce the auxiliary maps $\hat{f}_k:\mathbb{R}^n\times\mathbb{R}_{\geq 0}\to\mathbb{R}^n$, for $k\in\{1,2\}$, defined as
\begin{align}\label{constructionhatfk}
    \hat{f}_k(x,\tau)&:= \varphi(x) f_k(x,\tau), 
\end{align}
\end{subequations}
as well as the transformation $\Psi$, defined as
\begin{subequations}\label{eq:near_identity_transformation}
\begin{align}\label{nearidentity1}
    \pi_{1}\circ\Psi(x,\tau) &= \Phi(x,\tau), & \pi_{2}\circ\Psi(x,\tau)&= \tau,
\end{align}
for all $(x,\tau)\in\mathbb{R}^{n}\times\mathbb{R}_{\geq 0}$, where the map $\Phi$ is defined as follows:
\begin{align}\label{nearidentity2}
    \Phi(x,\tau) := x - \varepsilon\,{v}_1(x,\tau) - \varepsilon^2 {v}_2(x,\tau),
\end{align}
for all $(x,\tau)\in\mathbb{R}^{n}\times\mathbb{R}_{\geq 0}$, with
\begin{align}
    {v}_1(x,\tau)&:= \int_0^\tau \hat{f}_1(x,s)\,\text{d}s, \\
    {v}_2(x,\tau)&:= w(x,\tau) -  \D x v_1(x,\tau) v_1(x,\tau),\\
    w(x,\tau)&:= \int_0^\tau\left(\hat{f}_2(x,s) + \D x \hat{f}_1 (x,s)v_1(x,s)-\bar{f}(x)\right) \text{d}s,
\end{align}
and where the \emph{second-order average mapping} $\bar{f}$ is given by
\begin{align}\label{eq:averaged_vector_field}
    \bar{f}(x)&:= \frac{1}{2T}\int_0^T\left(2\hat{f}_2(x,\tau)+[v_1,\hat{f}_1](x,\tau)\right)\,\text{d}\tau,
\end{align}
for all $(x,\tau)\in\mathbb{R}^{n}\times\mathbb{R}_{\geq 0}$, with $T\in\mathbb{R}_{>0}$ being the same constant from \thref{asmp:vector_fields_condition_time_varying}, and $[v_1,\hat{f}_1]$ denoting the Lie bracket between the vector $v_1$ and $\hat{f}_1$, i.e.,
\begin{align*}
    [v_1,\hat{f}_1](x,\tau)&= \D{x}\hat{f}_1(x,\tau)v_1(x,\tau)-\D{x}v_1(x,\tau)\hat{f}_1(x,\tau).
\end{align*}
\end{subequations}
\begin{rem}\normalfont
The map $\Psi$ defined via \eqref{eq:near_identity_transformation} is an example of a (second-order) \textit{near-identity} transformation \cite{sanders2007averaging}, which is a standard tool in the averaging literature. The nomenclature stems from the fact that when $\varepsilon =0$, the transformation \eqref{eq:near_identity_transformation} reduces to the identity map on its domain and, by choosing $0<\varepsilon\ll 1$ sufficiently small, the transformation \eqref{eq:near_identity_transformation} can be made arbitrarily close to the identity map on bounded subsets of its domain \cite[Lemma 2.8.3]{sanders2007averaging}. Note that $\Psi$ depends (smoothly) on $\varepsilon$, but we suppress this dependency in the notation for brevity. \QEDB 
\end{rem}
%
%

\begin{figure}[t!]
    \centering
    \includegraphics[width=\linewidth]{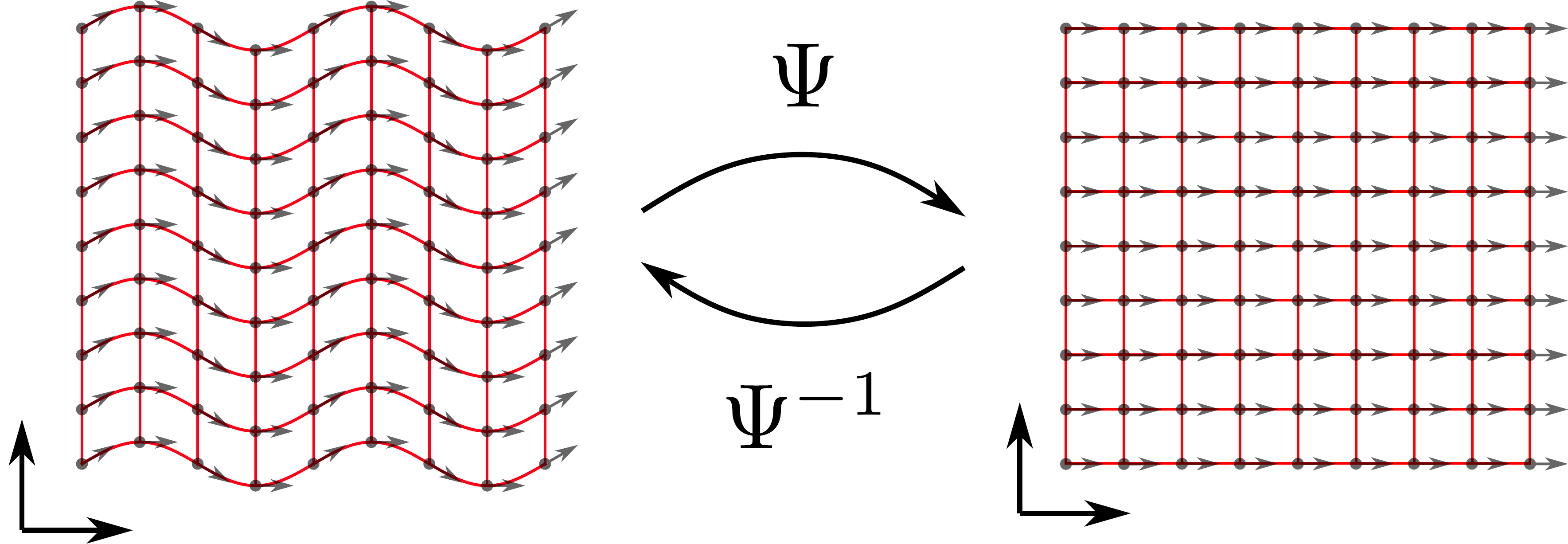}
    \caption{Visual depiction of the local effect of $\Psi$ on the integral curves of $f_{\varepsilon}$.}
    \label{fig:near_identity_transformation}
\end{figure}

Finally, for each $\epsilon>0$ and $\delta$ of the form \eqref{eq:deltas}, we also consider the closed sets
\begin{align}\label{eq:concentric_subsets}
    \mathcal{M}_j&:=\{x\in\mathbb{R}^n:\,|x|\geq \delta_j\}, & j&\in\{1,2,3\},
\end{align}
which satisfy $\mathcal{M}_1\supseteq\mathcal{M}_2\supseteq\mathcal{M}_3$. In fact, by construction, the case $\mathcal{M}_1=\mathcal{M}_2=\mathbb{R}^n$ can only occur if $\delta_1=0$.  We illustrate the function $\varphi$ and the sets $\mathcal{M}_j$, for $j\in\{1,2,3\}$, in Figure \ref{fig:bump_function}, and the local effect of the map $\Psi$, which will act as a coordinate transformation, in Figure \ref{fig:near_identity_transformation}.

\vspace{0.1cm}
The following proposition, key for our results, characterizes some useful properties of the function $\Psi$  and the pushforward under $\Psi$ of the vector field \eqref{eq:orig_system_timevaring}, c.f., Def. \ref{defn:diffeomorphism_action_on_ODE}.
%
%
%

%
\begin{prop}\thlabel{prop:near_identity_transformation}
     Suppose that \thref{asmp:vector_fields_condition_time_varying} holds, and let $\epsilon>0$ and $\delta$ satisfy \eqref{eq:deltas}. 
     Then, there exists $\varepsilon_0 \in\mathbb{R}_{>0}$, $L_\Psi \in\mathbb{R}_{>0}$, $L_g \in\mathbb{R}_{>0}$, and a $\mathcal{C}^0$ map $g:\mathbb{R}^n\times\mathbb{R}_{\geq0}\times[0,\varepsilon_0]\to\mathbb{R}^n$, such that for all $ \varepsilon\in(0,\varepsilon_0]$ the following holds:
     \vspace{0.1cm}
     \begin{enumerate}[(a),topsep=0pt,itemsep=-1ex,partopsep=1ex,parsep=1ex]
         \item\label{prop:near_identity_transformation_1} The map $\Psi:\mathbb{R}^{n}\times\mathbb{R}_{\geq 0}\to\mathbb{R}^{n}\times\mathbb{R}_{\geq 0}$ is a $\mathcal{C}^1$ diffeomorphism.
         \item\label{prop:near_identity_transformation_3}The map $\Psi$ and its inverse $\Psi^{-1}$ satisfy:
         \begin{align*}
                &\left|\pi_1\circ\Psi(0,\tau)\right|\leq L_\Psi \varepsilon,\\
                &\left|\pi_1\circ\Psi^{-1}(0,\tau)\right|\leq L_\Psi \varepsilon,\\
             &\left|\Psi(x_1,\tau)-\Psi(x_2,\tau)\right| \leq (1+L_\Psi \varepsilon)|x_1-x_2|,\\
             &\left|\Psi^{-1}(x_1,\tau)-\Psi^{-1}(x_2,\tau)\right| \leq (1+L_\Psi \varepsilon)|x_1-x_2|,
         \end{align*}
          for all $x_1,x_2\in\mathbb{R}^{n}$, and for all $\tau\in\mathbb{R}_{\geq 0}$.
        \item\label{prop:near_identity_transformation_4} For all $(x,\tau)\in\mathcal{M}_3\times\mathbb{R}_{\geq 0}$, we have $\Psi^{-1}(x,\tau)\in\mathcal{M}_2\times\mathbb{R}_{\geq 0}$.
        \item\label{prop:near_identity_transformation_5} The map $\Psi_*f_\varepsilon$ satisfies
        \begin{align}\label{eq:pushforward_map}
            \Psi_*f_\varepsilon(x,\tau)&= \bar{f}(x) + \varepsilon\,g(x,\tau,\varepsilon),
        \end{align}
        for all $(x,\tau)\in\mathcal{M}_3\times\mathbb{R}_{\geq 0}$. 
        \item\label{prop:near_identity_transformation_6} The map $g$ satisfies
        \begin{align}\label{eqnpertg}
            |g(x,\tau,\varepsilon)|&\leq L_g (|{x}|+1),
        \end{align}
        for all $({x},\tau,\varepsilon)\in\mathbb{R}^n\times\mathbb{R}_{\geq 0}\times [0,\varepsilon_0]$. \QEDB
     \end{enumerate}
\end{prop}

\noindent 
\textbf{Proof:} See Section \ref{proofProposition1}.

\begin{rem}\normalfont
Apart from the suitable smoothness and boundedness properties of $\Psi$,
\thref{prop:near_identity_transformation} asserts that system \eqref{eq:perturbed_avged_system_timevaring} can be seen as a perturbation of the nominal second-order average system
\begin{align}\label{eq:nominal_avged_sys}
    \dot{\bar{x}}&= \bar{f}(\bar{x}),~~~~~~\bar{x}\in\mathbb{R}^n,
\end{align}
for all $(x,\tau,\varepsilon)\in\mathcal{M}_3\times\mathbb{R}_{\geq 0}\times(0,\varepsilon_0]$, where $\bar{f}$ is given by \eqref{eq:averaged_vector_field}. By using this relationship, as well as the properties of $\Psi$, we can inform the stability analysis of system \eqref{eq:orig_system_timevaring} based on the stability properties of the nominal averaged system \eqref{eq:nominal_avged_sys}.\QEDB 
\end{rem}
\subsection{Global Stability via Second-Order Averaging}
To study the stability properties of \eqref{eq:orig_system_timevaring} via averaging, we make the following assumption on the average map $\bar{f}$.
%
%
\begin{asmp}\thlabel{asmp:Lyapunov_condition_timevarying}\normalfont
There exists  $\epsilon>0$, a vector $\delta$ satisfying \eqref{eq:deltas} with the same $\delta_1$ generated by \thref{asmp:vector_fields_condition_time_varying}, a $\mathcal{C}^1$-function $V:\mathbb{R}^n\to\mathbb{R}_{\geq0}$, $\alpha_i\in\mathcal{K}_\infty$, $c_i>0$, for $i\in\{1,2\}$, and a positive definite function $\phi:\mathbb{R}^n\to\mathbb{R}_{\geq0}$, such that the following holds:
\begin{enumerate}[(a)]
    \item \label{asmp:Lyapunov_condition_timevarying_1} For all $x\in\mathbb{R}^n$, we have that \begin{subequations}
        \begin{gather}
            \alpha_1(|x|) \leq V(x)\leq \alpha_2(|x|),\\
            |\nabla V(x)|\leq c_2\phi(x).
        \end{gather}
    \item \label{asmp:Lyapunov_condition_timevarying_2} For all $x\in\mathcal{M}_3$, we have that
    \begin{align}
            \left\langle\nabla V(x), \bar{f}(x)\right\rangle\leq -c_1\phi(x)^2.
    \end{align}
    \item \label{asmp:perturbation_norm_upper_bound} At least one of the following statements holds:
    \begin{enumerate}[(i),topsep=0pt,itemsep=-1ex,partopsep=1ex,parsep=1ex]
        \item \label{asmp:perturbation_norm_upper_bound_1} There exists $\bar{L}_{g}>0$, such that  $$|g(x,\tau,\varepsilon)|\leq \bar{L}_g (\phi(x) + 1),$$
        for all $(x,\tau,\varepsilon)\in\mathcal{M}_3\times\mathbb{R}_{\geq 0}\times[0,\varepsilon_0]$, where $g$ is the map generated by \thref{prop:near_identity_transformation}.
        \item\label{asmp:perturbation_norm_upper_bound_2} There exists $\alpha_3\in\mathcal{K}$, such that $\alpha_3(|x|)|x|\leq \phi(x)$.
    \end{enumerate}
   \end{subequations}
\end{enumerate}
 
 \QEDB
\end{asmp}
\begin{figure*}[ht]
    \centering
    \includegraphics[width=\linewidth]{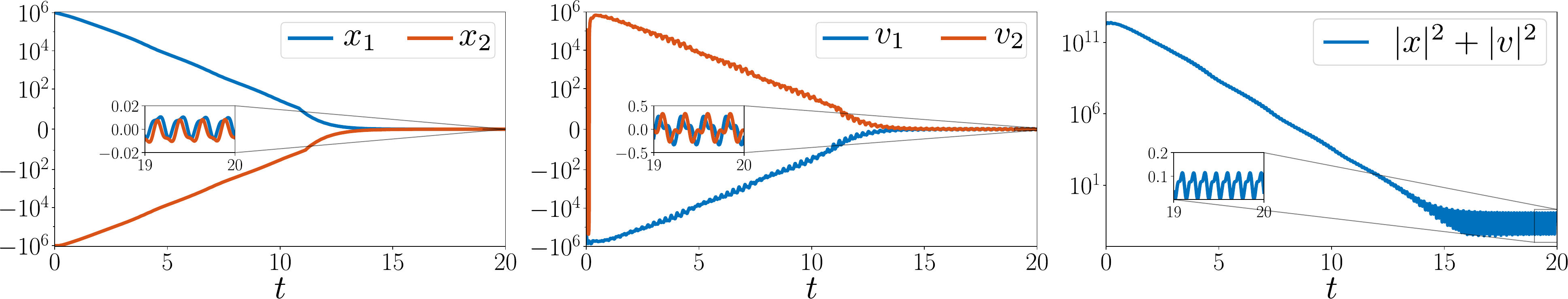}
    \caption[]{Numerical results for Example 1. The initial conditions are: $x(0)=(10^{6},-10^{6})$, $v(0)=(10^{3},-10^{3})$. The parameters are: $\varepsilon=1/\sqrt{8\pi}$, $\gamma_1=\frac{3}{4}$, $\gamma_2=1$. }
    \label{fig:clf_minimum_seeking_example_results}
\end{figure*}
The quadratic-type Lyapunov conditions in items \ref{asmp:Lyapunov_condition_timevarying_1} and \ref{asmp:Lyapunov_condition_timevarying_2} of \thref{asmp:Lyapunov_condition_timevarying} are identical to those studied in the literature of perturbed ODEs \cite[Section 9.1]{khalil2002nonlinear}. They imply that the origin is UGUB for the nominal average system \eqref{eq:nominal_avged_sys} \cite[Thm. 4.18]{khalil2002nonlinear}. However, without further restrictions on the rate of growth of the norm of the map $g$ in \eqref{eq:pushforward_map} relative to the map $\phi$, the perturbation $\varepsilon\,g$ may dominate the average map $\bar{f}$ far from the origin for any non-zero $\varepsilon$, thereby destroying global stability properties. To preclude this possibility, we impose the additional assumption in item \ref{asmp:perturbation_norm_upper_bound} of \thref{asmp:Lyapunov_condition_timevarying}.
%
%
\begin{rem}\normalfont
    Item \ref{asmp:perturbation_norm_upper_bound}-\ref{asmp:perturbation_norm_upper_bound_1} in \thref{asmp:Lyapunov_condition_timevarying} is automatically satisfied whenever the map $g$ is uniformly bounded. As we show in the proof of \thref{thm:es_bounded_gradient}, ES systems with bounded vector fields (see, e.g. \cite{scheinker2013model}) satisfy this condition under appropriate assumptions on the cost function. However, item \ref{asmp:perturbation_norm_upper_bound_1} leaves room for unbounded growth of the map $g$, provided that it can be dominated by the positive definite function $\phi$, for all $(x,\tau,\varepsilon)\in\mathcal{M}_3\times\mathbb{R}_{\geq 0}\times[0,\varepsilon_0]$. For example, item \ref{asmp:perturbation_norm_upper_bound}-\ref{asmp:perturbation_norm_upper_bound_1} in \thref{asmp:Lyapunov_condition_timevarying} automatically holds for the case $\phi(x)=|x|$ thanks to item \ref{prop:near_identity_transformation_6} in \thref{prop:near_identity_transformation}.
\end{rem}
\begin{rem}\normalfont
    Item \ref{asmp:perturbation_norm_upper_bound}-\ref{asmp:perturbation_norm_upper_bound_2} in \thref{asmp:Lyapunov_condition_timevarying} is automatically satisfied for the case $\phi(x)=|x|$. However, since $\alpha_3$ is an arbitrary $\mathcal{K}$ function, item \ref{asmp:perturbation_norm_upper_bound}-\ref{asmp:perturbation_norm_upper_bound_2} is a substantial relaxation of the local behavior of the function $\phi$ on any compact neighborhood of the origin.
\end{rem}
By leveraging the previous constructions and Proposition \ref{prop:near_identity_transformation}, we can now state the first main result of the paper. All the proofs are presented in Section \ref{sec:proofs}.
\begin{thm}\thlabel{thm:GUES_implies_GUPES_timevarying}
Suppose that Assumptions \ref{asmp:vector_fields_condition_time_varying}-\ref{asmp:perturbation_norm_upper_bound} hold.
    Then, for each $\epsilon>0$ and $\delta$ satisfying  \eqref{eq:deltas}, 
    there exists $\Delta>0$ such that the origing $x^*=0$ is $\Delta$-UGUB for system \eqref{pushforward}.
\end{thm}

We now provide several useful corollaries of \thref{thm:GUES_implies_GUPES_timevarying}. The first corollary concerns the stability properties of the original system \eqref{eq:orig_system_timevaring}.
\begin{cor}\thlabel{cor:GUES_implies_GUPES_timevarying_1}
    Suppose the assumptions of \thref{thm:GUES_implies_GUPES_timevarying} hold. Then, for each $\epsilon>0$ and $\delta$ satisfying \eqref{eq:deltas}, there exists $\tilde{\Delta}>0$ such that $x^*=0$ is $\tilde{\Delta}$-UGUB for system \eqref{eq:orig_system_timevaring}.
\end{cor}
\begin{cor}\thlabel{cor:GUES_implies_GUPES_timevarying_2}
    {Suppose that \thref{asmp:vector_fields_condition_time_varying} is satisfied for each $\delta_1>0$, and that there exists $V$ and $\phi$ satisfying the conditions of \thref{asmp:Lyapunov_condition_timevarying} for all $\epsilon>0$ and $\delta$ satisfying \eqref{eq:deltas}. Then, $x^*=0$ is UGpAS for system \eqref{eq:orig_system_timevaring}}. 
\end{cor}
\begin{cor}\thlabel{cor:GUES_implies_GUPES_timevarying_3}
    {Suppose that \thref{asmp:vector_fields_condition_time_varying} is satisfied for $\delta_1=0$, and that there exists a function $V$ and a function $\phi$ such that the conditions of  \thref{asmp:Lyapunov_condition_timevarying} hold for all $\epsilon>0$ and $\delta$ satisfying \eqref{eq:deltas}. Then, $x^*=0$ is UGpAS for system \eqref{eq:orig_system_timevaring}}.
\end{cor}

\begin{rem}\normalfont
     Corollary \ref{cor:GUES_implies_GUPES_timevarying_2} considers the situation in which \thref{asmp:vector_fields_condition_time_varying} is satisfied for each $\delta_1>0$ but might be violated for $\delta_1=0$. Such a situation arises when the vector fields defining system \eqref{eq:orig_system_timevaring} satisfy \thref{asmp:vector_fields_condition_time_varying} on any \emph{closed} subset of the set $\mathbb{R}^n\backslash\{0\}\times\mathbb{R}_{\geq 0}$, but strictly violate \thref{asmp:vector_fields_condition_time_varying} on $\mathbb{R}^n\times\mathbb{R}_{\geq 0}$. We illustrate this situation in \thref{exmp:non_lipschitz_vectorfield} below. \QEDB
\end{rem}
\begin{rem}\normalfont
Corollary \ref{cor:GUES_implies_GUPES_timevarying_3} considers the situation in which \thref{asmp:vector_fields_condition_time_varying} is satisfied for $\delta_1=0$, i.e. for all $(x,\tau)\in\mathbb{R}^n\times\mathbb{R}_{\geq 0}$, which differs significantly from the case considered in \thref{cor:GUES_implies_GUPES_timevarying_2}. \QEDB
\end{rem}
%
%
We conclude this section with a few illustrative examples. The first example serves the purpose of illustrating Corollary \ref{cor:GUES_implies_GUPES_timevarying_2}.
\begin{exmp}\thlabel{exmp:non_lipschitz_vectorfield}\normalfont
    Let $x\in\mathbb{R}$, $\tau\in\mathbb{R}_{\geq 0}$, and consider the dynamical system
    \begin{align}\label{eq:non_lipschitz_vectorfield}
        \dot{x}&= -|x|^{\frac{1}{2}}\text{sign}(x)\sin(\tau)^2, & \dot{\tau}&=\varepsilon^{-2},
    \end{align}
    which fits the structure of the class of systems \eqref{eq:orig_system_timevaring} with $f_1=0$ and $f_2(x,\tau)=-|x|^{\frac{1}{2}}\text{sign}(x)\sin(\tau)^2$. For any fixed $\delta_1>0$, there exists a constant $L_{\delta_1}>0$ such that, for all $x_1,x_2\in\mathbb{R}\backslash(-\delta_1,\delta_1)$ and all $\tau\in\mathbb{R}_{\geq 0}$, the function $f_2$ satisfies $|f_2(x_1,\tau)-f_2(x_2,\tau)|\leq L_{\delta_1}|x_1-x_2|.$ However, the constant $L_{\delta_1}>0$ tends to $+\infty$ in the limit $\delta_1\rightarrow 0$. In other words, there is no constant $L_0>0$ such that, for all $x_1,x_2\in\mathbb{R}$ and all $\tau\in\mathbb{R}_{\geq 0}$, we have $|f_2(x_1,\tau)-f_2(x_2,\tau)|\leq L_0|x_1-x_2|.$ Nevertheless, system \eqref{eq:non_lipschitz_vectorfield} satisfies \thref{asmp:vector_fields_condition_time_varying} for any $\delta_1>0$. Using formula \eqref{eq:averaged_vector_field}, we obtain that, for any choice of $\delta_1>0$, $\epsilon>0$, and $\delta$ satisfying \eqref{eq:deltas}, the corresponding nominal averaged system is given by
    \begin{align}
        \dot{\bar{x}}&= -\frac{1}{2}\varphi(\bar{x})|\bar{x}|^{\frac{1}{2}}\text{sign}(\bar{x}),
    \end{align}
    where $\varphi$ is the function defined in \eqref{bumpfunction}. Consider the function $V(x)=|x|^{\frac{3}{2}}$, which is $\mathcal{C}^1$, and the function $\phi(x)=|x|^{\frac{1}{2}}$, which is positive definite. Observe that the functions $V$ and $\phi$ satisfy item \ref{asmp:Lyapunov_condition_timevarying_1} in \thref{asmp:Lyapunov_condition_timevarying} with $c_1=\frac{3}{2}$. Direct computation gives 
    \begin{align}
        \nabla V(x)&=\frac{3}{2}|x|^{\frac{1}{2}}\,\text{sign}(x), & \left\langle \nabla V(x),\bar{f}(x)\right\rangle &= -\frac{3}{4}\varphi(x)\phi(x)^2.
    \end{align}
    Moreover, observe that, by construction, for any choice of $\delta_1>0$, $\epsilon>0$, and $\delta$ satisfying \eqref{eq:deltas}, we have that $\varphi(x)=1$, for all $|x|\geq \delta_2$. Consequently, it follows that the functions $V$ and $\phi$ satisfy item \ref{asmp:Lyapunov_condition_timevarying_2} in \thref{asmp:Lyapunov_condition_timevarying} with $c_2=\frac{3}{4}$, for any choice of $\delta_1>0$, $\epsilon>0$, and $\delta$ satisfying \eqref{eq:deltas}. Finally, notice that $f_1(x,\tau)=0$, and $|f_2(x,\tau)|\leq \phi(x)$, for all $x\in\mathbb{R}$ and all $\tau\in\mathbb{R}_{\geq 0}$. As a consequence, it can be shown that the function $g$ generated by \thref{prop:near_identity_transformation} is such that the function $\phi$ satisfies \ref{asmp:perturbation_norm_upper_bound}-\ref{asmp:perturbation_norm_upper_bound_1}, for any choice of $\delta_1>0$, $\epsilon>0$, and $\delta$ satisfying \eqref{eq:deltas}. Therefore, we have argued that system \eqref{eq:non_lipschitz_vectorfield} satisfies the assumptions of  \thref{cor:GUES_implies_GUPES_timevarying_2}. By invoking \thref{cor:GUES_implies_GUPES_timevarying_2}, we conclude that the point $x^*=0$ is UGpAS for system \eqref{eq:non_lipschitz_vectorfield}.
\end{exmp}
The second example, adapted from \cite{scheinker2012minimum}, shows that \emph{global} stabilization is possible via vibrational feedback control in certain systems with unknown control directions.
\begin{exmp}\normalfont
    Let $x=(x_1,x_2)\in\mathbb{R}^2,v=(v_1,v_2)\in\mathbb{R}^2$, $B\in\mathbb{R}^{2\times 2}$ such that $\text{rank}(B)=2$, and consider the dynamical system
    \begin{align}\label{eq:min_clf_exmp_1}
        \dot{x}= v,~~~~~~\dot{v}= B\hat{u},  
    \end{align}
    where $\hat{u}=(\hat{u}_1,\hat{u}_2)\in\mathbb{R}^2$ is the control input. The goal is to stabilize the equilibrium position $x=v=0$ for system \eqref{eq:min_clf_exmp_1} under the assumption that $B$ is unknown. To tackle this problem, we consider a model-free controller inspired by the ES systems studied in \cite{suttner2017exponential}. Namely, we let $\varepsilon\in\mathbb{R}_{>0}$, $\tau\in\mathbb{R}_{\geq 0}$, and consider the feedback law:
    \begin{subequations}\label{eq:min_clf_exmp_2}
        \begin{align}
            \hat{u}_1 &= \varepsilon^{-1} u_1(x,v,\tau), & \hat{u}_2 &= \varepsilon^{-1} u_2(x,v,\tau), & \dot{\tau}&= \varepsilon^{-2},
        \end{align}
        where the functions $u_i$ are given by
        \begin{align}
            u_1(x,v,\tau)&= \sqrt{2V(x,v)}\cos(\ln(V(x,v))+\tau), \\
           u_2(x,v,\tau)&=  \sqrt{4V(x,v)}\cos(\ln(V(x,v))+2\tau), 
        \end{align}
    \end{subequations}
    and where the function $V$ is taken as
    \begin{align}\label{eq:min_clf_exmp_3}
        V(x,v)&= |\gamma_1 x +\gamma_2v|^2+\frac{1}{2},
    \end{align}
    and the positive gains $\gamma_1$ and $\gamma_2$ are tuning parameters. It can be shown that the closed loop system defined by \eqref{eq:min_clf_exmp_1}-\eqref{eq:min_clf_exmp_3} satisfies \thref{asmp:vector_fields_condition_time_varying} for $\delta_1=0$ (see the proof of \thref{thm:es_unbounded_gradient} in Section \ref{sectionprooftheorem2}). Hence, we are allowed to pick $\delta_1=\delta_2=0$ and an arbitrary $\delta_3>0$. Using the formula \eqref{eq:averaged_vector_field}, we obtain that the nominal averaged system is given, for all $(x,v)\in\mathbb{R}^2\times\mathbb{R}^2$, by
    \begin{align}\label{eq:min_clf_exmp_4}
        \begin{pmatrix}
            \dot{\bar{x}} \\ \dot{\bar{v}}
        \end{pmatrix}&= A \begin{pmatrix}
            \bar{x} \\ \bar{v}
        \end{pmatrix} = \begin{pmatrix}
            0 & \text{I}\\-\gamma_1\gamma_2 BB^\top & -\gamma_2^2 BB^\top
        \end{pmatrix}\begin{pmatrix}
            \bar{x} \\ \bar{v}
        \end{pmatrix},
    \end{align}
    which turns out to be linear and time-invariant. If the matrix $A$ in \eqref{eq:min_clf_exmp_4} is Hurwitz, then system \eqref{eq:min_clf_exmp_4} is UGAS \cite[Theorem 4.5]{khalil2002nonlinear} and, by converse Lyapunov theorems \cite[Theorem 4.14]{khalil2002nonlinear}, it also satisfies \thref{asmp:Lyapunov_condition_timevarying} with $\phi(x,v)=|(x,v)|$, and $\alpha_3(r)=\tanh(r)$. Consequently, by invoking \thref{cor:GUES_implies_GUPES_timevarying_3} we conclude that the closed-loop system defined by \eqref{eq:min_clf_exmp_1}-\eqref{eq:min_clf_exmp_3} renders the origin UGpAS. Figure \ref{fig:clf_minimum_seeking_example_results} shows the behavior exhibited by the trajectories of the system. In all the simulations, we used $B = (1,1;1,-1)$. We remark that, although we take $B$ here as a constant matrix, a similar result can be established when $B$ is time-varying under suitable uniform persistence of excitation conditions, see \cite{scheinker2012minimum}. 
    \QEDB 
\end{exmp}
\section{Applications to Extremum Seeking Systems}\label{sec:es_section}
\label{sec:applicationsES}
In this section, we leverage the averaging results established in Theorem \ref{thm:GUES_implies_GUPES_timevarying} and Corollaries \ref{cor:GUES_implies_GUPES_timevarying_1}-\ref{cor:GUES_implies_GUPES_timevarying_3} to study uniform global practical asymptotic stability (UGpAS) for a class of ES systems of the form \eqref{eq:es_system_open_loop}.
\subsection{Main Assumptions}
To guarantee that the (open-loop) amplitudes of the exploration signals in \eqref{eq:es_system_open_loop} have access to all directions in the parameter space, we consider the following assumption on the vectors $b_{i,j}$.

\begin{asmp}\thlabel{asmp:vector_fields_conditions_es}\normalfont
        There exists $\gamma>0$, such that the vectors $b_{i,j}$ satisfy $\sum_{i=1}^r\sum_{j=1}^2 \left(b_{i,j}^\top v\right)^2 \geq \gamma |v|^2$, for all $v\in\mathbb{R}^n$. \QEDB 
\end{asmp}
    We also make the following regularity assumption on the cost functions $J$ and the drift term $b_0$. In all cases, we assume that  $J^\star:=\inf_{x\in\mathbb{R}^n} J(x)>-\infty$, and that $J^\star=J(x^\star)$ for some unique $x^\star\in\mathbb{R}^n$. Similar conditions were used in \cite{grushkovskaya2018class,suttner2017exponential} to analyze ES systems with (local) asymptotic stability properties.
    \begin{asmp}\thlabel{asmp:drift_term_es}\normalfont
        The following holds:
        \vspace{0.2cm}
        \begin{enumerate}[(a),topsep=0pt,itemsep=-0.5ex,partopsep=1ex,parsep=1ex]
    \item $J(x)>J(x^\star)$, for all $x\neq x^\star$.
        \item $\nabla J(x)=0$ if and only if $x=x^\star$.
        \item There exists $L_J>0$ such that $|\nabla^2 J(x)|\leq L_J$.
        \item There exists $\kappa_3>0$ such that $|b_0(x)|\leq \kappa_3|\nabla J(x)|$, for all $x\in\mathbb{R}^n$.
        \item There exists $L_0>0$ such that $|b_0(x_1)-b_0(x_2)|\leq L_0|x_1-x_2|$,  for all $x_1,x_2\in\mathbb{R}^n$. 
        \end{enumerate}
        \vspace{-0.3cm}
        \QEDB 
    \end{asmp}
      \begin{rem}\normalfont
    Items (a)-(b) in \thref{asmp:cost_function_conditions_es_1} are standard in ES problems \cite{tan2006non,krstic2000stability}. Similarly, item (c) is equivalent to the assumption that  $\nabla J$ is $L_J$-globally Lipschitz \cite[Lemma 1.2.2]{NesterovsBook}, which is satisfied by quadratic maps, which are typically studied in ES problems \cite{ariyur2003real}. Finally, note that items (d)-(e) are relevant only when the drift term $b_0$ in \eqref{eq:es_system_open_loop} is not zero. However, in most ES systems this term is set to zero. \QEDB
    \end{rem}
    Next, we consider two classes of cost functions $J:\mathbb{R}^n\to\mathbb{R}$ that we seek to globally minimize via the dynamics \eqref{eq:es_system_open_loop}. 
\begin{asmp}\thlabel{asmp:cost_function_conditions_es_1}\normalfont
        The cost $J$ is a radially unbounded $\mathcal{C}^2$-function and there exists $\alpha\in\mathcal{K}$ such that at least one of the following statements holds:
        \vspace{0.2cm}
        \begin{enumerate} [a),topsep=0pt,itemsep=-0.5ex,partopsep=1ex,parsep=1ex]

        \item For all $x\in\mathbb{R}^n$, we have  $$\alpha_J(|x-x^\star|)^2|x-x^\star|^2\leq|\nabla J(x)|^2.$$
        \item There exists $M_J>0$ such that  $$\alpha_J(|x-x^\star|)^2\leq|\nabla J(x)|^2\leq M_J^2,$$ for all $x\in\mathbb{R}^n$. \QEDB
        \end{enumerate}
        %
    \end{asmp}
    \begin{rem}\label{remnonconvex}\normalfont
    As shown in \thref{lem:strong_convexity_implies_assumption} in the Appendix, item (a) in \thref{asmp:cost_function_conditions_es_1} is satisfied by any strongly convex $\mathcal{C}^2$-function with a globally Lipschitz gradient. This family of functions includes quadratic cost functions having a positive definite Hessian, which are common in ES. However, as shown in the next example, convexity of the cost function $J$ is not needed to satisfy \thref{asmp:cost_function_conditions_es_1}.  \QEDB
    \end{rem}
    %
    %
    \begin{exmp}\thlabel{exmp:es_1_illustration}\normalfont
    Let $n=2$, $x^\star=(10^{10},-10^{10})$, and let the cost function $J:\mathbb{R}^2\rightarrow\mathbb{R}$ be given by
    \begin{equation}\label{nonconvexcost}
    J(x):=|x-x^\star|^2+3\sin\left(|x-x^\star|\right)^2+1,
    \end{equation}
    which is not convex \cite[pp.4]{karimi2016linear}. However, as shown in Lemma \ref{lem:exmp_2_cost_function} in the Appendix, $J$ satisfies items (a)-(c) in \thref{asmp:drift_term_es} and item (a) in \thref{asmp:cost_function_conditions_es_1} with $L_J=20$ and class-$\mathcal{K}$ function $\alpha_J(s)=0.5\tanh(s)$.  \QEDB 
    \end{exmp}
    %

            \begin{figure*}[ht]
        \centering
        \includegraphics[width=\linewidth]{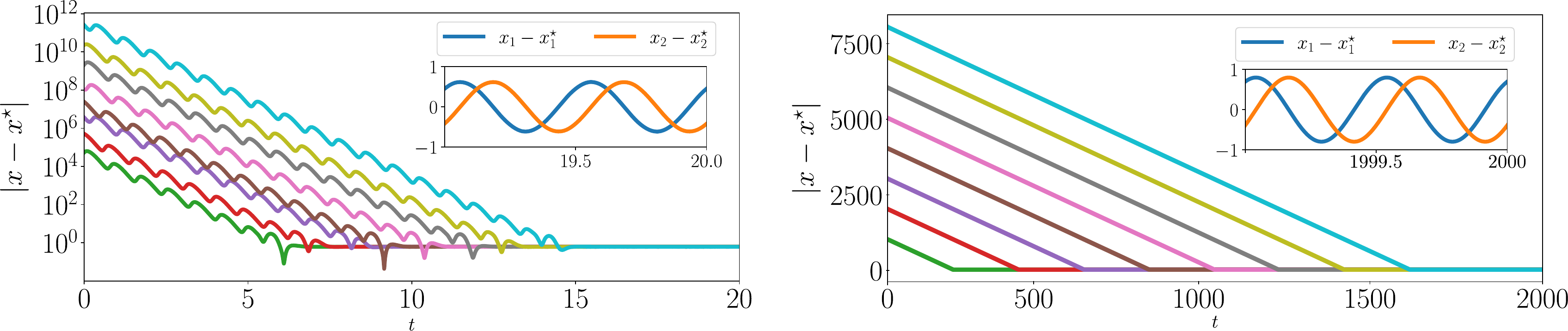}
        \caption[]{Numerical results for Example 4 (left) and Example 5 (right). The insets in the top right of the figures depict the quasi-steady state in the vicinity of $x^\star$.}
        \label{fig:es_example_results}
    \end{figure*}
    
    The following example considers a cost function $J$ obtained as a regularization of the vector norm function, which satisfies item b) in \thref{asmp:cost_function_conditions_es_1}.
    \begin{exmp}\thlabel{exmp:es_2_illustration}\normalfont
    Let $n=2$, $x^\star=(10^3,-10^{3})$, and let the cost function $J:\mathbb{R}^n\rightarrow\mathbb{R}$ be given by
    \begin{equation}\label{seconcost}
    J(x):=|x-x^\star|\tanh(|x-x^\star|) -100.
    \end{equation}
    Therefore, it can be directly verified that the function $J$ satisfies items (a)-(c) in \thref{asmp:drift_term_es} and item b) in \thref{asmp:cost_function_conditions_es_1} with $L_J=3$, $M_J=2$, and the class-$\mathcal{K}$ function $\alpha_J(s)=\tanh(s)$.\QEDB 
    \end{exmp}
    \subsection{ES Dynamics with Linear Growth}
    We now consider two different algorithms of the form \eqref{eq:es_system_open_loop} that are able to achieve global ES. The first algorithm that we consider, initially introduced in \cite{suttner2017exponential}, can be written as system \eqref{eq:es_system_open_loop} with the following time-varying feedback law:
    \begin{subequations}\label{eq:es_law_1}
    \begin{align}
        u_{i,1}(J(x),\tau)&:= \begin{cases}
            \sqrt{2\omega_i J(x)}\cos(\ln(J(x))+\omega_i\tau) & J(x)> 0\\
            0 & J(x)\leq 0
            \
        \end{cases},\\
        u_{i,2}(J(x),\tau)&:= \begin{cases}
            \sqrt{2\omega_i J(x)}\sin(\ln(J(x))+\omega_i\tau) & J(x)>0\\
            0 & J(x)\leq 0
        \end{cases},
    \end{align}
    \end{subequations}
    where $\omega_i\in\mathbb{Q}_{>0}$, such that $\omega_i\neq\omega_j$ for $i\neq j$. 
    
    Using the coordinate shift $\tilde{x}=x-x^\star$, system \eqref{eq:es_system_open_loop} can be written as follows:
    \begin{subequations}\label{eq:closed_loop_es_orig_sys}
        \begin{align}
            \dot{\tilde{x}}&= b_0(\tilde{x}+x^\star) + \varepsilon^{-1}\sum_{\substack{i=1}}^{r}\sum_{j=1}^2 b_{i,j}u_{i,j}(J(\tilde{x}+x^\star),\tau), \\ \dot{\tau}&= \varepsilon^{-2},
        \end{align}
    \end{subequations}
    which has the same form as \eqref{eq:orig_system_timevaring}. 

   \vspace{0.2cm}
    The following theorem is the second main result of this paper.
    \begin{thm}\thlabel{thm:es_unbounded_gradient}
        Suppose that \thref{asmp:vector_fields_conditions_es,asmp:drift_term_es} hold with $\gamma > \kappa_3$. Then, if item a) in \thref{asmp:cost_function_conditions_es_1} holds:
        \begin{enumerate}[(a)]
            \item The origin $\tilde{x}=0$ is UGUB for system \eqref{eq:closed_loop_es_orig_sys} under the feedback law \eqref{eq:es_law_1}.
            \item If $J^\star\geq 0$, then the origin $\tilde{x}=0$ is UGpAS for system \eqref{eq:closed_loop_es_orig_sys} under the feedback law \eqref{eq:es_law_1}. 
        \end{enumerate}
        \vspace{-0.6cm}
        \QEDB
    \end{thm}
    The novelty of Theorem \ref{thm:es_unbounded_gradient} is to establish \emph{global} bounds of the form \eqref{KLbound1001}-\eqref{UUbound} for the ES dynamics \eqref{eq:es_system_open_loop} with feedback law \eqref{eq:es_law_1}. As discussed in Example \ref{exmp:es_1_illustration}, such bounds can be obtained even when $J$ is not convex.
    
\begin{exmp}\thlabel{exmp:es_1}\normalfont
    (Example \ref{exmp:es_1_illustration} continued) Let $\omega_1=1$, $r=1$, and consider the ES system \eqref{eq:es_system_open_loop} with destabilizing drift $b_0(x)=\frac{1}{2}(x-x^\star)$, and constant vectors $b_{1,1}=(1,0)$, $b_{1,2}=(0,1)$. Notice that in this case $\gamma = 1$, and it can be shown that $\kappa_3=0.8$. Therefore, item d) in \thref{asmp:drift_term_es} is satisfied. The feedback law is given by \eqref{eq:es_law_1}, with cost function \eqref{nonconvexcost}. Since all the assumptions of Theorem \ref{thm:es_unbounded_gradient} are satisfied and $J^{\star}>0$, we conclude that $x^{\star}$ is UGpAS. Numerical simulation results are shown in Fig. \ref{fig:es_example_results}. In the figure, we present simulations obtained from various randomly generated initial conditions and using $\varepsilon=1/\sqrt{4\pi}$. As observed, all trajectories converge to a neighborhood of the optimizer. \QEDB 
\end{exmp}
    \subsection{ES Dynamics with Bounded Control}
    The second ES algorithm that we consider can also be written as system \eqref{eq:es_system_open_loop} with feedback law:
    \begin{subequations}\label{eq:es_law_2}
    \begin{align}
        u_{i,1}(J(x),\tau)&:= \sqrt{2\omega_i}\cos(J(x)+\omega_i\tau),\\
        u_{i,2}(J(x),\tau)&:= \sqrt{2\omega_i}\sin(J(x)+\omega_i\tau).
    \end{align}
    \end{subequations}
    The semi-global practical stability properties of these systems have been studied in \cite{scheinker2013model,scheinker2014extremum}. These algorithms are characterized by uniformly bounded vector fields, which are suitable for applications with actuator constraints.  Note that using the coordinate shift $\tilde{x}=x-x^\star$, system \eqref{eq:es_system_open_loop} with the feedback law \eqref{eq:es_law_2} can also be written as \eqref{eq:closed_loop_es_orig_sys}.

    \vspace{0.1cm}
    The following theorem is the third main result of this paper. 
    %
    %
    \begin{thm}\thlabel{thm:es_bounded_gradient}
        Suppose that \thref{asmp:vector_fields_conditions_es,asmp:drift_term_es} hold. If $\gamma > \kappa_3$ and item b) in \thref{asmp:cost_function_conditions_es_1} is satisfied, then the origin $\tilde{x}=0$ is UGpAS for system \eqref{eq:closed_loop_es_orig_sys} under the feedback law \eqref{eq:es_law_2}. \QEDB
    \end{thm}
The novelty of Theorem \ref{thm:es_bounded_gradient} compared to the results of \cite{scheinker2013model,scheinker2014extremum}, is to establish a global bound of the form \eqref{KLbound1001} for all solutions of the system, albeit under stronger assumptions on the cost functions.

We conclude this section by presenting a numerical example that illustrates the application of Theorem  \ref{thm:es_bounded_gradient}.

%
\begin{exmp}\thlabel{exmp:es_2}\normalfont
    (Example \ref{exmp:es_2_illustration} continued) Let $\omega_1=1$, $r=1$, and consider the ES system \eqref{eq:es_system_open_loop} with $b_0(x)=(0,0)$, $b_{1,1}=(2,0)$, $b_{1,2}=(0,2)$, and the feedback law \eqref{eq:es_law_2}. Notice that in this case $\kappa_3=0$ and $\gamma = 2$, so the assumption that $\kappa_3<\gamma$ holds trivially. We consider the cost function \eqref{seconcost}, which satisfies the required Assumptions to apply Theorem \ref{thm:es_bounded_gradient}. We simulate the system from randomly generated initial conditions with $\varepsilon=1/\sqrt{4\pi}$. Numerical simulation results are shown in Fig. \ref{fig:es_example_results}. Since item b) in \thref{asmp:cost_function_conditions_es_1} restricts the gradient to be uniformly bounded, the convergence rate that emerges is slower compared to the convergence rate of the ES dynamics of \thref{exmp:es_1}. However, \thref{thm:es_bounded_gradient} still asserts UGpAS of the minimizer. \QEDB 
\end{exmp} 
\section{Proofs}\label{sec:proofs}
In this section, we present the proofs of the main results.
\subsection{Proof of \thref{prop:near_identity_transformation}}
\label{proofProposition1}
    For the sake of clarity, we divide the proof into several key lemmas. The following lemma is a direct consequence of \thref{asmp:vector_fields_condition_time_varying} and the properties of the map $\varphi$ summarized in \thref{lem:hole_function_properties} (see also \cite[Lemmas 2.20-2.22]{lee2012smooth}).
    \begin{lem}\thlabel{lem:vector_fields_condition_time_varying} Let the assumptions of \thref{prop:near_identity_transformation} be satisfied. Then, for all $k\in\{1,2\}$, the following holds:
        \begin{enumerate}[(a),topsep=0pt,itemsep=-1ex,partopsep=1ex,parsep=1ex]
            \item\label{lem:vector_fields_condition_time_varying_1} The map $\hat{f}_k$ is $\mathcal{C}^0$ and there exist positive constants $\hat{L}_{k}$ such that for all $x_1,x_2\in\mathbb{R}^n$ and all $\tau\in\mathbb{R}_{\geq 0}$, the map $\hat{f}_k$ satisfies
            \begin{gather*}
                \left|\hat{f}_k(x_1,\tau)-\hat{f}_k(x_2,\tau)\right|\leq \hat{L}_k |x_1-x_2|.
            \end{gather*}
            \item \label{lem:vector_fields_condition_time_varying_2} For all $(x,\tau)\in\mathbb{R}^n\times\mathbb{R}_{\geq 0}$, and for the same $T\in\mathbb{R}_{>0}$ from item \ref{asmp:vector_fields_condition_time_varying_2} in \thref{asmp:vector_fields_condition_time_varying}, we have that
            \begin{align*}
                \hat{f}_k(x,\tau+T)&=\hat{f}_k(x,\tau), &  {\textstyle\int}_{0}^T \hat{f}_1(x,\tau) d\tau = 0. 
            \end{align*}
            \item\label{lem:vector_fields_condition_time_varying_3} The map $\hat{f}_k$ is $\mathcal{C}^{3-k}$ with respect to $x$ on $\mathbb{R}^{n}\times\mathbb{R}_{\geq 0}$.
            \item\label{lem:vector_fields_condition_time_varying_4} There exists a positive constant $\hat{L}_{3} $ such that, for all $(x_1,\tau_1),(x_2,\tau_2)\in\mathbb{R}^n\times\mathbb{R}_{\geq 0}$, the map $\D{x} \hat{f}_1\cdot\hat{f}_1$, satisfies
            \begin{align*}
                |\D{x} \hat{f}_1(x_1,\tau_1)\cdot \hat{f}_1(x_1,\tau_2)-\D{x} &\hat{f}_1(x_2,\tau_1)\cdot \hat{f}_1(x_2,\tau_2)|\\
                &\leq \hat{L}_{3} |x_1-x_2|.
            \end{align*}
            \item\label{lem:vector_fields_condition_time_varying_5} For all $(x,\tau)\in\mathcal{M}_2$, we have that $\hat{f}_k(x,\tau)= f_k(x,\tau).$
            \item\label{lem:vector_fields_condition_time_varying_6} For all $(x,\tau)\in\text{cl}\left(\mathbb{R}^n\backslash\mathcal{M}_1\right)\times\mathbb{R}_{\geq 0}$, we have that $\hat{f}_k(x,\tau)= 0.$
        \end{enumerate}
    \end{lem}
    \noindent\textbf{Proof: } If $\delta_1=\delta_2=0$, then $\varphi(x)=1$ and all the properties follow directly by \eqref{constructionhatfk} and Assumption \ref{asmp:vector_fields_condition_time_varying}. On the other hand, suppose that $\delta_2 > \delta_1$. We begin by proving items \ref {lem:vector_fields_condition_time_varying_5} and \ref{lem:vector_fields_condition_time_varying_6}. By \thref{lem:hole_function_properties}, $\varphi(x)=1$, for all $x\in\mathcal{M}_2$. Hence, by definition, the map $\hat{f}_k$ satisfies $\hat{f}_k(x,\tau)=f_k(x,\tau)$, for all $(x,\tau)\in\mathcal{M}_2\times\mathbb{R}_{\geq 0}$, which proves item \ref {lem:vector_fields_condition_time_varying_5}. Similarly, by \thref{lem:hole_function_properties}, $\varphi(x)=0$, for all $x\in\delta_1\mathbb{B}$. Hence, by construction, the map $\hat{f}_k$ also satisfies $\hat{f}_k(x,\tau)=0$, for all $(x,\tau)\in\delta_1\mathbb{B}\times\mathbb{R}_{\geq 0}$, which proves item \ref{lem:vector_fields_condition_time_varying_6}.

    Next, we prove item \ref{lem:vector_fields_condition_time_varying_1}. By \thref{lem:hole_function_properties}, the map $\varphi(x)$ is $\mathcal{C}^\infty$. Hence, by item \ref{asmp:vector_fields_condition_time_varying_1} in \thref{asmp:vector_fields_condition_time_varying}, the definition of the map $\hat{f}_k$ implies that $\hat{f}_k(\cdot,\cdot)$ is $\mathcal{C}^0$ on $\mathbb{R}^n\times\mathbb{R}_{\geq 0}$. Since $\hat{f}_k(x,\tau)=f_k(x,\tau)$, for all $(x,\tau)\in\mathcal{M}_2\times\mathbb{R}_{\geq 0}$, it follows that $\hat{f}_k$ inherits all the properties of $f_k$ in the domain $\mathcal{M}_2\times\mathbb{R}_{\geq 0}$. In particular, items \ref{asmp:vector_fields_condition_time_varying_1} and \ref{asmp:vector_fields_condition_time_varying_3} in \thref{asmp:vector_fields_condition_time_varying} imply that $\D{x}f_k$ is well-defined and satisfies the bound $|\D{x}f_k(x,\tau)|\leq L_k$, for all $(x,\tau)\in\mathcal{M}_2\times\mathbb{R}_{\geq 0}$. Consequently, $\D{x}\hat{f}_k(x,\tau)$ also satisfies the bound $|\D{x}\hat{f}_k(x,\tau)|\leq L_k$, for all $(x,\tau)\in\mathcal{M}_2\times\mathbb{R}_{\geq 0}$. Similarly, since $\hat{f}_k(x,\tau)=0$, for all $(x,\tau)\in\delta_1\mathbb{B}\times\mathbb{R}_{\geq 0}$, it follows that $\D{x}\hat{f}_k(x,\tau)=0$ is well-defined and satisfies $\D{x}\hat{f}_k(x,\tau)=0$, for all $(x,\tau)\in\delta_1\mathbb{B}\times\mathbb{R}_{\geq 0}$. Finally, the definition of the map $\hat{f}_k$ implies that, for all $(x,\tau)\in\text{cl}\left(\mathcal{M}_1\backslash\mathcal{M}_2\right)\times\mathbb{R}_{\geq 0}$, we have
    \begin{align*}
        \D{x}\hat{f}_k(x,\tau)&= \varphi(x)\D{x}f_k(x,\tau) + f_k(x,\tau)\nabla \varphi(x)^\top,
    \end{align*}
    which is continuous. Since $\text{cl}\left(\mathcal{M}_1\backslash\mathcal{M}_2\right)$ is compact and, for all $x\in\mathbb{R}^n$, $\D{x}\hat{f}_k(x,\cdot):\mathbb{R}_{\geq 0}\rightarrow \mathbb{R}^{n\times n}$ is $\mathcal{C}^0$ and periodic, it follows that there exists a constant $\tilde{L}_k\in\mathbb{R}_{>0}$ such that, for all $(x,\tau)\in\text{cl}\left(\mathcal{M}_1\backslash\mathcal{M}_2\right)\times\mathbb{R}_{\geq 0}$, the Jacobian $\D{x}\hat{f}_k$ satisfies the upper bound $|\D{x}\hat{f}_k(x,\tau)|\leq \tilde{L}_k$. Let $\hat{L}_k:=\max\{L_k,\tilde{L}_k\}$. Then, for all $(x,\tau)\in\mathbb{R}^n\times\mathbb{R}_{\geq 0}$, the Jacobian $\D{x}\hat{f}_k$ is well-defined and satisfies the upper bound $|\D{x}\hat{f}_k(x,\tau)|\leq \hat{L}_k$. Consequently, for all $x_1,x_2\in\mathbb{R}^n$, and $\forall \tau\in\mathbb{R}_{\geq 0}$, the map $\hat{f}_k$ satisfies
    \begin{align*}
        |\hat{f}_k(x_1,\tau)-\hat{f}_k(x_2,\tau)|&\leq \hat{L}_k|x_1-x_2|,
    \end{align*}
    which proves item \ref{lem:vector_fields_condition_time_varying_1}. Item \ref{lem:vector_fields_condition_time_varying_2} follows directly from the definition of the map $\hat{f}_k$.

    Next, we prove item \ref{lem:vector_fields_condition_time_varying_3}. Since the map $\varphi$ is $\mathcal{C}^\infty$, the definition of the map $\hat{f}_k$ implies that it inherits all the smoothness properties of $f_k$ in the domain $\mathcal{M}_1\times\mathbb{R}_{\geq 0}$. In particular, item \ref{asmp:vector_fields_condition_time_varying_3} in \thref{asmp:vector_fields_condition_time_varying} implies that $\hat{f}_k(\cdot,\tau)$ is $\mathcal{C}^{3-k}$ on the closed set $\mathcal{M}_1$, for all $\tau\in\mathbb{R}_{\geq 0}$. On the other hand, since $\hat{f}_k(x,\tau)=0$ for all $(x,\tau)\in\delta_1\mathbb{B}\times\mathbb{R}_{\geq 0}$, it follows that $\hat{f}_k(\cdot,\tau)$ is $\mathcal{C}^{\infty}$ on the open set $\mathbb{R}^n\backslash\mathcal{M}_1$, for all $\tau\in\mathbb{R}_{\geq 0}$. It follows that $\hat{f}_k(\cdot,\tau)$ is $\mathcal{C}^{3-k}$ on $\mathbb{R}^n$, for all $\tau\in\mathbb{R}_{\geq 0}$, which proves item \ref{lem:vector_fields_condition_time_varying_3}.
    
    Next, we prove item \ref{lem:vector_fields_condition_time_varying_4}. For $(x,\tau_1,\tau_2)\in\mathbb{R}^n\times\mathbb{R}_{\geq 0}\times\mathbb{R}_{\geq 0}$, define the maps
    \begin{align*}
        F(x,\tau_1,\tau_2)&:=\D{x}\hat{f}_1(x,\tau_1)f_1(x,\tau_2), \\
        \hat{F}(x,\tau_1,\tau_2)&:=\varphi(x)F(x,\tau_1,\tau_2) = \D{x}\hat{f}_1(x,\tau_1)\hat{f}_1(x,\tau_2).
    \end{align*}
    Since $\hat{f}_1(\cdot,\tau_1)$ is $\mathcal{C}^2$ on $\mathbb{R}^n$, for all $\tau_1\in\mathbb{R}_{\geq 0}$, and $f_1(\cdot,\tau_1)$ is $\mathcal{C}^2$ on $\mathcal{M}_1$, for all $\tau_1\in\mathbb{R}_{\geq 0}$, it follows that the map $F(\cdot,\tau_1,\tau_2)$ is $\mathcal{C}^{1}$ on $\mathcal{M}_1$, for all $(\tau_1,\tau_2)\in\mathbb{R}_{\geq 0}\times\mathbb{R}_{\geq 0}$. In addition, since the map $\varphi$ is $\mathcal{C}^\infty$, it follows that the map $F(\cdot,\tau_1,\tau_2)$ is also $\mathcal{C}^{1}$ on $\mathcal{M}_1$, for all $(\tau_1,\tau_2)\in\mathbb{R}_{\geq 0}\times\mathbb{R}_{\geq 0}$. From \thref{lem:hole_function_properties}, $\varphi(x)=1$, for all $x\in\mathcal{M}_2$. Hence, by definition, the map $\hat{F}$ satisfies $\hat{F}(x,\tau_1,\tau_2)=F(x,\tau_1,\tau_2)$, for all $(x,\tau_1,\tau_2)\in\mathcal{M}_2\times\mathbb{R}_{\geq 0}\times\mathbb{R}_{\geq 0}$, which means that $\hat{F}$ inherits all the properties of the map $F$ in the domain $\mathcal{M}_2\times\mathbb{R}_{\geq 0}\times\mathbb{R}_{\geq 0}$. In particular, from items \ref{asmp:vector_fields_condition_time_varying_3} and \ref{asmp:vector_fields_condition_time_varying_4} in \thref{asmp:vector_fields_condition_time_varying}, $\D{x}F$ is well-defined and satisfies the bound $|\D{x}F(x,\tau_1,\tau_2)|\leq L_3$, for all $(x,\tau_1,\tau_2)\in\mathcal{M}_2\times\mathbb{R}_{\geq 0}\times\mathbb{R}_{\geq 0}$, which implies that $\D{x}\hat{F}$ also satisfies the bound $|\D{x}\hat{F}(x,\tau_1,\tau_2)|\leq L_3$, for all $(x,\tau_1,\tau_2)\in\mathcal{M}_2\times\mathbb{R}_{\geq 0}\times\mathbb{R}_{\geq 0}$. From \thref{lem:hole_function_properties}, $\varphi(x)=0$, for all $x\in\delta_1\mathbb{B}$. Hence, by definition, the map $\hat{F}$ satisfies $\hat{F}(x,\tau_1,\tau_2)=0$, for all $(x,\tau_1,\tau_2)\in\delta_1\mathbb{B}\times\mathbb{R}_{\geq 0}\times\mathbb{R}_{\geq 0}$, which implies that $\hat{F}(\cdot,\tau_1,\tau_2)$ is $\mathcal{C}^{1}$ on $\mathbb{R}^n$, for all $(\tau_1,\tau_2)\in\mathbb{R}_{\geq 0}\times\mathbb{R}_{\geq 0}$, and that $\D{x}\hat{F}$ satisfies $\D{x}\hat{F}(x,\tau_1,\tau_2)=0$, for all $(x,\tau_1,\tau_2)\in\delta_1\mathbb{B}\times\mathbb{R}_{\geq 0}\times\mathbb{R}_{\geq 0}$. Finally, the definition of the map $\hat{F}$ implies that, for all $(x,\tau)\in\text{cl}\left(\mathcal{M}_1\backslash\mathcal{M}_2\right)\times\mathbb{R}_{\geq 0}$, we have
    \begin{align*}
        \D{x}\hat{F}(x,\tau_1,\tau_2)&= \varphi(x)\D{x}F(x,\tau_1,\tau_2) + F(x,\tau_1,\tau_2)\nabla \varphi(x)^\top,
    \end{align*}
    which is continuous. Since $\text{cl}\left(\mathcal{M}_1\backslash\mathcal{M}_2\right)$ is compact and, for all $(x,\tau_2)\in\mathbb{R}^n\times\mathbb{R}_{\geq 0}$, $\D{x}\hat{F}(x,\cdot,\tau_2):\mathbb{R}_{\geq 0}\rightarrow \mathbb{R}^{n\times n}$ is $\mathcal{C}^0$ and periodic, and for all $(x,\tau_1)\in\mathbb{R}^n\times\mathbb{R}_{\geq 0}$, $\D{x}\hat{F}(x,\tau_1,\cdot):\mathbb{R}_{\geq 0}\rightarrow \mathbb{R}^{n\times n}$ is $\mathcal{C}^0$ and periodic, it follows that there exists a constant $\tilde{L}_3\in\mathbb{R}_{>0}$ such that, for all $(x,\tau_1,\tau_2)\in\text{cl}\left(\mathcal{M}_1\backslash\mathcal{M}_2\right)\times\mathbb{R}_{\geq 0}\times\mathbb{R}_{\geq 0}$, the Jacobian $\D{x}\hat{F}$ satisfies the upper bound $|\D{x}\hat{F}(x,\tau_1,\tau_2)|\leq \tilde{L}_3$. Let $\hat{L}_3:=\max\{L_3,\tilde{L}_3\}$. Then, for all $(x,\tau_1,\tau_2)\in\mathbb{R}^n\times\mathbb{R}_{\geq 0}\times\mathbb{R}_{\geq 0}$, the Jacobian $\D{x}\hat{F}$ satisfies the upper bound $|\D{x}\hat{F}(x,\tau_1,\tau_2)|\leq \hat{L}_3$. Consequently, for all $x_1,x_2\in\mathbb{R}^n$, and $\forall (\tau_1,\tau_2)\in\mathbb{R}_{\geq 0}\times\mathbb{R}_{\geq 0}$, the map $\hat{F}$ satisfies
    \begin{align*}
        |\hat{F}(x_1,\tau_1,\tau_2)-\hat{F}(x_2,\tau_1,\tau_2)|&\leq \hat{L}_3|x_1-x_2|,
    \end{align*}
    which proves item \ref{lem:vector_fields_condition_time_varying_4} and concludes the proof of the Lemma. \hfill$\blacksquare$
    \begin{lem}\label{lem:Lipschitz_conditions}
        Let the assumptions of \thref{prop:near_identity_transformation} be satisfied. Then, the following holds: 
        \begin{enumerate}[(a),topsep=0pt,itemsep=-1ex,partopsep=1ex,parsep=1ex]
            \item \label{lem:Lipschitz_conditions_1} The maps $\bar{f}$ and $v_k$, for $k\in\{1,2\}$, are $\mathcal{C}^1$ on $\mathbb{R}^n\times\mathbb{R}_{\geq 0}$.
            \item \label{lem:Lipschitz_conditions_2} There exist $\bar{L},L_{v,k}>0$, for $k\in\{1,2\}$, such that:
            \begin{align*}
                &|\bar{f}(x_1)-\bar{f}(x_2)|\leq \bar{L} |x_1-x_2|,\\
                &|v_k(x_1,\tau)-v_k(x_2,\tau)|\leq L_{v,k} |x_1-x_2|,~~\forall~k\in\{1,2\},
            \end{align*}
            for all $x_1,x_2\in\mathbb{R}^n$ and all $ \tau\in\mathbb{R}_{\geq 0}$.
        \end{enumerate}
    \end{lem}
    \noindent\textbf{Proof: }
        First, we prove item a). Since $v_1$ is the integral of $\hat{f}_1$ with respect to $\tau$, it follows from item \ref{lem:vector_fields_condition_time_varying_3} in \thref{lem:vector_fields_condition_time_varying} that the map $v_1$ is $\mathcal{C}^2$ in $x$, and $\mathcal{C}^1$ in $\tau$. In addition, since $\bar{f}$ is a multiple of the definite integral with respect to $\tau$ of the terms $\D x v_1\,\hat{f}_1$, $\D x \hat{f}_1\,v_1$, and $\hat{f}_2$, and, from item \ref{lem:vector_fields_condition_time_varying_3} in \thref{lem:vector_fields_condition_time_varying}, all those of terms are $\mathcal{C}^1$ in $x$, it follows that $\bar{f}$ is $\mathcal{C}^1$. Moreover, since $v_2$ is the sum of the term $\D x v_1 v_1$, which is $\mathcal{C}^1$ in all arguments, and the integral with respect to $\tau$ of the terms $\D x v_1\,\hat{f}_1$, $\D x \hat{f}_1\,v_1$, and $\hat{f}_2$, which are all, from item \ref{lem:vector_fields_condition_time_varying_3} in \thref{lem:vector_fields_condition_time_varying}, $\mathcal{C}^1$ in $x$ and $\mathcal{C}^0$ in $\tau$, it follows that $v_2$ is $\mathcal{C}^1$ in all arguments.
        
        Next, we prove item b). From the definition of the map $v_1$ and item \ref{lem:vector_fields_condition_time_varying_1} in \thref{lem:vector_fields_condition_time_varying}, we have that
        \begin{align*}
            |v_1(x_1,\tau)-v_1(x_2,\tau)|&=\left|{\textstyle\int}_0^{\tau}\left(\hat{f}_1(x_1,s)-\hat{f}_1(x_2,s)\right) \text{d}s\right|\\
            &\leq \hat{L}_{1} \tau|x_1-x_2|\leq T\hat{L}_1  |x_1-x_2|,
        \end{align*}
        for all $x_1,x_2\in\mathbb{R}^n$ and all $\tau\in[0,T]$. In addition, from item \ref{lem:vector_fields_condition_time_varying_2} in \thref{lem:vector_fields_condition_time_varying}, $v_1$ is periodic in $\tau$. It follows that, for all $x_1,x_2\in\mathbb{R}^n$ and all $\tau\in\mathbb{R}_{\geq 0}$, we have
        \begin{align*}
            |v_1(x_1,\tau)-v_1(x_2,\tau)|&\leq L_{v,1} |x_1-x_2|, & L_{v,1} &:=T\hat{L}_1 .
        \end{align*}
        From the definition of $v_1$, and by interchanging matrix multiplication with the integral, we have that
        \begin{align*}
            \D x \hat{f}_1(x,\tau)v_1(x,\tau)&=  {\textstyle\int}_0^\tau\D x \hat{f}_1(x,\tau) \hat{f}_1(x,s)\text{d}s
        \end{align*}
        From item \ref{lem:vector_fields_condition_time_varying_4} in \thref{lem:vector_fields_condition_time_varying}, we have that
        \begin{align*}
            |\D x \hat{f}_1(x_1&,\tau) v_1(x_1,\tau)-\D x \hat{f}_1(x_2,\tau) v_1(x_2,\tau)|\\
            &=
            \left|{\textstyle\int}_0^\tau \left(\D x \hat{f}_1(x_1,\tau) \hat{f}_1(x_1,s)-\D x \hat{f}_1(x_2,\tau) \hat{f}_1(x_2,s)\right)\text{d}s\right|\\
            &\leq \hat{L}_{3} \tau |x_1-x_2| \leq T\hat{L}_{3} |x_1-x_2|.
        \end{align*}
        for all $x_1,x_2\in\mathbb{R}^n$ and all $\tau\in[0,T]$. In addition, from item \ref{lem:vector_fields_condition_time_varying_2} in \thref{lem:vector_fields_condition_time_varying}, $v_1$ and $\hat{f}_1$ are periodic in $\tau$. It follows that, for all $x_1,x_2\in\mathbb{R}^n$ and all $\tau\in\mathbb{R}_{\geq 0}$, we have that
        \begin{align*}
            |\D x \hat{f}_1(x_1&,\tau) v_1(x_1,\tau)-\D x \hat{f}_1(x_2,\tau) v_1(x_2,\tau)|\leq T\hat{L}_{3} |x_1-x_2|.
        \end{align*}
        From the definition of $v_1$, using Leibniz's rule, and by interchanging matrix multiplication with the integral, we have that
        \begin{align*}
            \D x v_1(x,\tau)\hat{f}_1(x,\tau)&=  {\textstyle\int}_0^\tau\D x \hat{f}_1(x,s) \hat{f}_1(x,\tau)\text{d}s
        \end{align*}
        From item \ref{lem:vector_fields_condition_time_varying_4} in \thref{lem:vector_fields_condition_time_varying}, we have that
        \begin{align*}
            |\D x v_1(x_1&,\tau) \hat{f}_1(x_1,\tau)-\D x v_1(x_2,\tau) \hat{f}_1(x_2,\tau)|\\
            &=
            \left|{\textstyle\int}_0^\tau \left(\D x \hat{f}_1(x_1,s) \hat{f}_1(x_1,\tau)-\D x \hat{f}_1(x_2,s) \hat{f}_1(x_2,\tau)\right)\text{d}s\right|\\
            &\leq \hat{L}_{3} \tau |x_1-x_2| \leq T\hat{L}_{3} |x_1-x_2|.
        \end{align*}
        for all $x_1,x_2\in\mathbb{R}^n$ and all $\tau\in[0,T]$. In addition, from item \ref{lem:vector_fields_condition_time_varying_2} in \thref{lem:vector_fields_condition_time_varying}, $v_1$ and $\hat{f}_1$ are periodic in $\tau$. It follows that, for all $x_1,x_2\in\mathbb{R}^n$ and all $\tau\in\mathbb{R}_{\geq 0}$, we have that
        \begin{align*}
            |\D x v_1(x_1&,\tau) \hat{f}_1(x_1,\tau)-\D x v_1(x_2,\tau) \hat{f}_1(x_2,\tau)|\leq T\hat{L}_{3} |x_1-x_2|.
        \end{align*}
        From the definition of $v_1$, using Leibniz's rule, and interchanging matrix multiplication with the integral, we have that
        \begin{align*}
            \D x v_1(x,\tau)v_1(x,\tau)&=  {\textstyle\int}_0^\tau{\textstyle\int}_0^\tau\D x \hat{f}_1(x,s) \hat{f}_1(x,\sigma)\text{d}s\,\text{d}\sigma
        \end{align*}
        From item \ref{lem:vector_fields_condition_time_varying_4} in \thref{lem:vector_fields_condition_time_varying}, we have that
        \begin{align*}
            |&\D x v_1(x_1,\tau) v_1(x_1,\tau)-\D x v_1(x_2,\tau) v_1(x_2,\tau)|\\
            &=
            \left|{\textstyle\int}_0^\tau {\textstyle\int}_0^\tau\left(\D x \hat{f}_1(x_1,s) \hat{f}_1(x_1,\sigma)-\D x \hat{f}_1(x_2,s) \hat{f}_1(x_2,\sigma)\right)\text{d}s\,\text{d}\sigma\right|\\
            &\leq \hat{L}_{3} \tau^2 |x_1-x_2| \leq T^2\hat{L}_{3} |x_1-x_2|.
        \end{align*}
        for all $x_1,x_2\in\mathbb{R}^n$ and all $\tau\in[0,T]$. In addition, from item \ref{lem:vector_fields_condition_time_varying_2} in \thref{lem:vector_fields_condition_time_varying}, $v_1$ is periodic in $\tau$. It follows that, for all $x_1,x_2\in\mathbb{R}^n$ and all $\tau\in\mathbb{R}_{\geq 0}$, we have that
        \begin{align*}
            |\D x v_1(x_1&,\tau) v_1(x_1,\tau)-\D x v_1(x_2,\tau) v_1(x_2,\tau)|\leq T^2\hat{L}_{3} |x_1-x_2|.
        \end{align*}
        Finally, note that, for all $x_1,x_2\in\mathbb{R}^n$, we have that
        \begin{align*}
            T|\bar{f}(x_1)-&\bar{f}(x_2)|\leq {\textstyle\int}_0^T|\hat{f}_2(x_1,\tau)-\hat{f}_2(x_2,\tau)| \text{d}\tau\\
            &+{\textstyle\int}_0^T|\D x v_1(x_1,\tau) \hat{f}_1(x_1,\tau)-\D x v_1(x_2,\tau) \hat{f}_1(x_2,\tau)|\text{d}\tau \\
            &+{\textstyle\int}_0^T|\D x \hat{f}_1(x_1,\tau) v_1(x_1,\tau)-\D x \hat{f}_1(x_2,\tau) v_1(x_2,\tau)| \text{d}\tau\\
            &\leq T^2\left(\hat{L}_2+2\hat{L}_3\right)|x_1-x_2|,
        \end{align*}
        and, for all $x_1,x_2\in\mathbb{R}^n$ and all $\tau\in\mathbb{R}_{\geq 0}$, we have that
        \begin{align*}
            |v_2(x_1,\tau)-&v_2(x_2,\tau)|\leq {\textstyle\int}_0^\tau|\hat{f}_2(x_1,\tau)-\hat{f}_2(x_2,\tau)| \text{d}\tau\\
            &+{\textstyle\int}_0^\tau|\bar{f}(x_1)-\bar{f}(x_2)| \text{d}\tau\\
            &+{\textstyle\int}_0^T|\D x \hat{f}_1(x_1,\tau) v_1(x_1,\tau)-\D x \hat{f}_1(x_2,\tau) v_1(x_2,\tau)| \text{d}\tau\\
            &+|\D x v_1(x_1,\tau) v_1(x_1,\tau)-\D x v_1(x_2,\tau) v_1(x_2,\tau)|\\
            &\leq T^2\left(2\hat{L}_2+4\hat{L}_3\right)|x_1-x_2|,
        \end{align*}
        The proof of the Lemma is concluded by defining $\bar{L}:= T(\hat{L}_2+2\hat{L}_3)$ and $L_{v,2}:=T^2(2\hat{L}_2+4\hat{L}_3)$. \hfill$\blacksquare$
        
        \begin{lem}\label{lem:near_identity_transformation_properties} Let the assumptions of \thref{prop:near_identity_transformation} be satisfied. Then, there exists $\varepsilon_0,L_\Psi>0$ such that for all $\varepsilon\in[0,\varepsilon_0]$, the following holds:
        \begin{enumerate}[(a),topsep=0pt,itemsep=-1ex,partopsep=1ex,parsep=1ex]
            \item \label{lem:near_identity_transformation_properties_1} $\Psi$ is a $\mathcal{C}^1$-diffeomorphism on $\mathbb{R}^n\times\mathbb{R}_{\geq 0}$.
            \item \label{lem:near_identity_transformation_properties_2} For all $\tau\in\mathbb{R}_{\geq 0}$, the diffeomorphism $\Psi$ and its inverse $\Psi^{-1}$ satisfy
            \begin{align*}
                \left|\pi_1\circ\Psi(0,\tau)\right|&\leq L_\Psi \varepsilon, & \left|\pi_1\circ\Psi^{-1}(0,\tau)\right|&\leq L_\Psi \varepsilon.
            \end{align*}
            \item \label{lem:near_identity_transformation_properties_3} For all $x_1,x_2\in\mathbb{R}^{n}$ and for all $\tau\in\mathbb{R}_{\geq 0}$, the map $\Psi$ and its inverse $\Psi^{-1}$ satisfy
            \begin{gather*}
                 \left|\Psi(x_1,\tau)-\Psi(x_2,\tau)\right| \leq (1+L_\Psi \varepsilon)|x_1-x_2|,\\
                 \left|\Psi^{-1}(x_1,\tau)-\Psi^{-1}(x_2,\tau)\right| \leq (1+L_\Psi \varepsilon)|x_1-x_2|.
            \end{gather*}
            \item \label{lem:near_identity_transformation_properties_4} For all $(x,\tau)\in\mathcal{M}_3\times\mathbb{R}_{\geq 0}$, $\Psi^{-1}(x,\tau)\in\mathcal{M}_2$.
        \end{enumerate}
    \end{lem}
    \noindent\textbf{Proof:} We first establish item a). Let $(x,\tau)$ and $(\tilde{x},\tilde{\tau})$ be any two points in $\mathbb{R}^n\times\mathbb{R}_{\geq 0}$ and suppose that $\Psi(x,\tau)=\Psi(\tilde{x},\tilde{\tau})$. Then, by construction, we have $\tau=\tilde{\tau}$, and
        \begin{align*}
            |x-\tilde{x}|&\leq  {\textstyle\sum}_{i=1}^2\varepsilon^i|{v}_i(x,\tau)-{v}_i(\tilde{x},\tau))|.
        \end{align*}
        From item \ref{lem:Lipschitz_conditions_2} in Lemma \ref{lem:Lipschitz_conditions}, we obtain that
        \begin{align*}
            |x-\tilde{x}|&\leq \varepsilon\,(L_{v,1}  + L_{v,2} \varepsilon)|x-\tilde{x}|
        \end{align*}
        Let $\bar{\varepsilon}_1:=\min\{1,1/(2(L_{v,1}  + L_{v,2} ))$ and note that for all $\varepsilon\in[0,\bar{\varepsilon}_1 ]$ we have that $|x-\tilde{x}|\leq \frac{1}{2}|x-\tilde{x}|$, which can only happen if $|x-\tilde{x}|=0$. Therefore, for all $(x,{\tau})\in\mathbb{R}^n\times\mathbb{R}_{\geq 0}$, for all $(\tilde{x},\tilde{\tau})\in\mathbb{R}^n\times\mathbb{R}_{\geq 0}$, and for all $\varepsilon\in[0,\bar{\varepsilon}_1 ]$, $\Psi(x,\tau)=\Psi(\tilde{x},\tilde{\tau})\implies x=\tilde{x},$ and $\tau=\tilde{\tau}$, which in turn implies that the map $\Psi$ is injective on $\mathbb{R}^n\times\mathbb{R}_{\geq 0}$. 
        Next, for each $(\tilde{x},\tau)\in\mathbb{R}^n\times\mathbb{R}_{\geq 0}$, define the map $\tilde{\Phi}:\mathbb{R}^n\rightarrow \mathbb{R}^n$: 
        \begin{align*}
            \tilde{\Phi}(x)&= \tilde{x} + x-\Phi(x,\tau).
        \end{align*}
        By direct computation
        \begin{align*}
            \tilde{\Phi}(x)&= \tilde{x} + {\textstyle\sum}_{k=1}\varepsilon^k v_k(x,\tau).
        \end{align*}
        Now let $x_1,x_2\in\mathbb{R}^n$ be any two points and observe that
        \begin{align*}
            \left|\tilde{\Phi}(x_1)-\tilde{\Phi}(x_2)\right|&\leq {\textstyle\sum}_{k=1}\varepsilon^k |v_k(x_1,\tau)-v_k(x_2,\tau)|.
        \end{align*}
        From item \ref{lem:Lipschitz_conditions_2} in Lemma \ref{lem:Lipschitz_conditions}, we obtain that
        \begin{align*}
            \left|\tilde{\Phi}(x_1)-\tilde{\Phi}(x_2)\right|&\leq \varepsilon\,(L_{v,1}  + L_{v,2} \varepsilon)|x_1-x_2|.
        \end{align*}
        Let $\varepsilon\in[0,\bar{\varepsilon}_1 ]$. Then, $\left|\tilde{\Phi}(x_1)-\tilde{\Phi}(x_2)\right|\leq \frac{1}{2}|x_1-x_2|$ for all $x_1,x_2\in\mathbb{R}^n$, which implies that $\tilde{\Phi}$ is a contraction. Thus,  $\tilde{\Phi}$ has a unique fixed point \cite[Lemma C.35]{lee2012smooth}, which implies that for all $\varepsilon\in[0,\bar{\varepsilon}_1 ]$, for all $\tilde{x}\in\mathbb{R}^n$ and for all $\tau\in\mathbb{R}_{\geq 0}$, there exists a unique point $x\in\mathbb{R}^n$ such that 
        \begin{equation}\label{importantonto}
        \tilde{x}=\pi_1\circ\Psi(x,\tau)=\Phi(x,\tau).
        \end{equation}
        In other words, $\Psi$ is onto, and therefore a bijection on $\mathbb{R}^n\times\mathbb{R}_{\geq 0}$. 
        
        From item \ref{lem:Lipschitz_conditions_1} in Lemma \ref{lem:Lipschitz_conditions}, we know that the map $\Psi$ is $\mathcal{C}^1$ on $\mathbb{R}^n\times\mathbb{R}_{\geq0}$, and so its Jacobian $\D{}\Psi$ is well-defined, and given by
        \begin{align*}
            \D{} \Psi&= \begin{pmatrix}\D x\Phi & \D \tau\Phi\\0 & 1\end{pmatrix},
        \end{align*}
        where the Jacobians $\D x\Phi$ and $ \D \tau\Phi$ are given by
        \begin{align*}
            {\D x\Phi}&= \text{I} -  {\textstyle\sum}_{k=1}^2\varepsilon^k\D x{v}_k, & 
            {\D \tau\Phi}&= -  {\textstyle\sum}_{k=1}^2\varepsilon^k\D \tau{v}_k.
        \end{align*}
         For each $i\in\{1,\dots,n\}$ and $k\in\{1,2\}$, let $R^i_{v,k}:\mathbb{R}^n\times\mathbb{R}_{\geq 0}\rightarrow \mathbb{R}_{\geq 0}$ be given by 
        \begin{align*}
            R^i_{v,k}(x,\tau)&= {\textstyle\sum}_{j=1,\,j\neq i}^n \left|\left(D_x v_k(x,\tau)\right)_{ij}\right|,
        \end{align*}
        where $\left(D_x v_k(x,\tau)\right)_{ij}$ are the entries of the matrix $D_x v_k(x,\tau)$. From items \ref{lem:Lipschitz_conditions_1} and \ref{lem:Lipschitz_conditions_2} in Lemma \ref{lem:Lipschitz_conditions}, the maps $v_k$ are $\mathcal{C}^1$ and globally Lipschitz with respect to $x$, uniformly in $\tau$. It follows that, for all $(x,\tau)\in\mathbb{R}^n\times\mathbb{R}_{\geq 0}$ and for $k\in\{1,2\}$, we have that $\left|D_x v_k(x,\tau)\right|\leq L_{v,k} $. In particular, for all $(x,\tau)\in\mathbb{R}^n\times\mathbb{R}_{\geq 0}$, for all $i,j\in\{1,\dots,n\}$, and for all $k\in\{1,2\}$, we have that $\left|\left(D_x v_k(x,\tau)\right)_{ij}\right|\leq L_{v,k} $. As such, for all $(x,\tau)\in\mathbb{R}^n\times\mathbb{R}_{\geq 0}$, for all $i\in\{1,\dots,n\}$, and for all $k\in\{1,2\}$, we have that $0\leq R^i_{v,k}(x,\tau)\leq (n-1)L_{v,k} $. 

        Next, note that the entries of $\D{x}\Phi$ are given by:
        \begin{align*}
            \left(\D{x}\Phi(x,\tau)\right)_{ii}&= 1-{\textstyle\sum}_{k=1}^2\varepsilon^k \left(\D{x}v_k(x,\tau)\right)_{ii}, & \forall i&\in\{1,\dots,n\} \\
            \left(\D{x}\Phi(x,\tau)\right)_{ij}&=-{\textstyle\sum}_{k=1}^2\varepsilon^k \left(\D{x}v_k(x,\tau)\right)_{ij}, & \forall i\neq j&\in\{1,\dots,n\},
        \end{align*}
        Consequently, for all $(x,\tau)\in \mathbb{R}^n\times\mathbb{R}_{\geq 0}$, for all $\varepsilon\in \mathbb{R}_{\geq 0}$, and for all 
        $i\in\{1,\dots,n\}$, we have
        \begin{align*}
            R^i(x,\tau)&:={\textstyle\sum}_{j\neq i=1}^n\left|\left(\D{x}\Phi(x,\tau)\right)_{ij}\right|\\
            &\leq  {\textstyle\sum}_{k=1}^2\varepsilon^k{\textstyle\sum}_{j\neq i=1}^n\left|\left(\D{x}v_k(x,\tau)\right)_{ij}\right|\\
            &=  {\textstyle\sum}_{k=1}^2\varepsilon^k R^i_{v,k}(x,\tau) \leq (n-1){\textstyle\sum}_{k=1}^2\varepsilon^kL_{v,k}. 
        \end{align*}
        Similarly, we have
        \begin{align*}
            1-{\textstyle\sum}_{k=1}^2\varepsilon^k L_{v,k} \leq \left(\D{x}\Phi(x,\tau)\right)_{ii}\leq 1 + {\textstyle\sum}_{k=1}^2\varepsilon^k L_{v,k} .
        \end{align*}
        Let $\bar{\varepsilon}_2 :=\min\{1, 1/(4 n (L_{v,1} +L_{v,2} ))\}$. Then, it follows that for all $\varepsilon\in[0,\bar{\varepsilon}_2 ]$, all $(x,\tau)\in\mathbb{R}^n\times\mathbb{R}_{\geq 0}$, and all $i\in\{1,\dots,n\}$:
        \begin{gather*}
            \frac{3}{4}\leq 1-\tilde{L}_{\Psi} \varepsilon\leq \left(\D{x}\Phi(x,\tau)\right)_{ii}\leq 1+\tilde{L}_{\Psi} \varepsilon\leq \frac{5}{4}, \\
            0\leq R^i(x,\tau)\leq \tilde{L}_{\Psi} \varepsilon\leq \frac{1}{4}.
        \end{gather*}
        where $\tilde{L}_{\Psi} :=n (L_{v,1} +L_{v,2} )$. 
        By applying the Ger{\v s}hgorin circle theorem \cite[p.269]{bernstein2009matrix}, we obtain that for all $\varepsilon\in[0,\bar{\varepsilon}_2 ]$, for all $(x,\tau)\in\mathbb{R}^n\times\mathbb{R}_{\geq 0}$, the eigenvalues of the Jacobian matrix $\D{x}\Phi(x,\tau)$ are contained in the compact interval $[1-2\tilde{L}_{\Psi} \varepsilon,1+2\tilde{L}_{\Psi} \varepsilon]\subset [1/2,3/2]$. 
        
        Then, we have the following claim, proved in  \ref{profclaim1}.
        \begin{claim}\label{refclaim1}
            For all $\varepsilon\in[0,\bar{\varepsilon}_2 ]$, for all $(x,\tau)\in\mathbb{R}^n\times\mathbb{R}_{\geq 0}$, the Jacobian matrix $\D{x}\Phi(x,\tau)$ is invertible and there exists a constant $L_\Psi \in\mathbb{R}_{>0}$ such that
        \begin{align*}
             \left|\D{x}\Phi(x,\tau)\right|&\leq 1+L_\Psi \varepsilon\leq 2, & \left|\D{x}\Phi(x,\tau)^{-1}\right|&\leq 1+L_\Psi \varepsilon\leq 2.
        \end{align*}
        and $D\Psi(x,\tau)^{-1}$ is well-defined and given by
        \begin{align}\label{inversetransfor}
            \D{} \Psi(x,\tau)^{-1}&= \begin{pmatrix}\D x\Phi(x,\tau)^{-1} & -\D x\Phi(x,\tau)^{-1}\D \tau\Phi(x,\tau)\\0 & 1\end{pmatrix},
        \end{align}
        for all $(x,\tau)\in\mathbb{R}^n\times\mathbb{R}_{\geq0}$. \QEDB
        \end{claim}
         Let $\bar{\varepsilon}:=\min\{\bar{\varepsilon}_1 ,\bar{\varepsilon}_2 \}$. For all $\varepsilon\in[0,\bar{\varepsilon} ]$, the map $\Psi$ is bijective and, for all $(x,\tau)\in\mathbb{R}^n\times\mathbb{R}_{\geq 0}$, the Jacobian $\D{}\Psi(x,\tau)$ is invertible and its inverse is continuous. Thus, by invoking the global rank theorem \cite[Theorem 4.14]{lee2012smooth}, we conclude that, for all $ \varepsilon\in[0,\bar{\varepsilon} ]$, the map $\Psi$ is a $\mathcal{C}^1$-diffeomorphism.

         \vspace{0.1cm}
         To establish item (b), we note that by \cite[Proposition C.4]{lee2012smooth} the Jacobian of $\Psi^{-1}$ is given by
        \begin{align*}
            \D{}\Psi^{-1}(x,\tau)&= \left(\D{}\Psi\circ\Psi^{-1}(x,\tau)\right)^{-1}.
        \end{align*}
        %
        
        \noindent From the definition of $\Psi$, for all $(\tau,\varepsilon)\in\mathbb{R}_{\geq 0}\times[0,\bar{\varepsilon}]$, we have that $\pi_1\circ\Psi(0,\tau)=-\sum_{k=1}^2\varepsilon^k v_k(0,\tau)$. Since $v_k(0,\tau)$ is periodic in $\tau$, there exists $M_k >0$ such that $|v_k(0,\tau)|\leq M_k $, for all $\tau$. Without loss of generality, we assume that $L_\Psi \geq M_1 +M_2 $. It follows that, for all $(\tau,\varepsilon)\in\mathbb{R}_{\geq 0}\times[0,\bar{\varepsilon}]$, we have $ |\pi_1\circ\Psi(0,\tau)|\leq L_\Psi \varepsilon$. Since $\Psi$ is a diffeomorphism, we have
        \begin{align*}
            0=\pi_1\circ\Psi\circ\Psi^{-1}(0,\tau)=\pi_1\circ\Psi^{-1}(0,\tau) -  {\textstyle\sum}_{k=1}^2\varepsilon^k{v}_k\circ\Psi^{-1}(0,\tau),
        \end{align*}
        which implies that $|\pi_1\circ\Psi^{-1}(0,\tau)|\leq|{\textstyle\sum}_{k=1}^2\varepsilon^k{v}_k\circ\Psi^{-1}(0,\tau)|$. Similar reasoning allows us to conclude that, without loss of generality, we may assume that $|\pi_1\circ\Psi^{-1}(0,\tau)|\leq L_\Psi \varepsilon$.

        \vspace{0.1cm}
        To establish item (c), note that, for all $x_1,x_2\in\mathbb{R}^n$, for all $\tau\in\mathbb{R}_{\geq 0}$, and for all $\varepsilon\in[0,\bar{\varepsilon} ]$, we obtain via Hadamard's Lemma \cite[Lemma 2.8]{nestruev2003smooth} that
        \begin{align*}
            \left|\Psi(x_1,\tau)-\Psi(x_2,\tau)\right|&\leq \left|\text{J}_{\Psi}(x_1,x_2,\tau)\right||x_1-x_2|,\\
            \left|\Psi^{-1}(x_1,\tau)-\Psi^{-1}(x_2,\tau)\right|&\leq \left|\text{J}_{\Psi^{-1}}(x_1,x_2,\tau)\right||x_1-x_2|,
        \end{align*}
        where 
        \begin{align*}
            \text{J}_{\Psi}(x_1,x_2,\tau)&:={\textstyle\int}_{0}^1\D{x}\Phi(x_2+\lambda(x_1-x_2),\tau) \text{d}\lambda \\
            \text{J}_{\Psi^{-1}}(x_1,x_2,\tau)&:={\textstyle\int}_{0}^1\D{x}\Phi^{-1}(x_2+\lambda(x_1-x_2),\tau) \text{d}\lambda,
        \end{align*}
        and where we used 
        \begin{align*}
            \D{x}\Phi^{-1}(x,\tau)&:=\left(\D{x}\Phi\circ\Psi^{-1}(x,\tau)\right)^{-1}.
        \end{align*}
        It follows that
        \begin{subequations}
            \label{eq:lipschitz_near_identity}
            \begin{align}
                \left|\Psi(x_1,\tau)-\Psi(x_2,\tau)\right|&\leq (1+L_\Psi \varepsilon)|x_1-x_2|, \\
                \left|\Psi^{-1}(x_1,\tau)-\Psi^{-1}(x_2,\tau)\right|&\leq (1+L_\Psi \varepsilon)|x_1-x_2|.
            \end{align}
        \end{subequations}
        for all $x_1,x_2\in\mathbb{R}^n$, for all $\tau\in\mathbb{R}_{\geq 0}$, and for all $\varepsilon\in[0,\bar{\varepsilon} ]$.

        Finally, to establish item (d), note that, since $\Psi$ is a diffeomorphism, we have
        \begin{align*}
            x=\pi_1\circ\Psi\circ\Psi^{-1}(x,\tau)=\pi_1\circ\Psi^{-1}(x,\tau) -  {\textstyle\sum}_{k=1}^2\varepsilon^k{v}_k\circ\Psi^{-1}(x,\tau).
        \end{align*}
        Therefore, for all $(x,\tau,\varepsilon)\in\mathbb{R}^n\times\mathbb{R}_{\geq 0}\times[0,\bar{\varepsilon} ]$:
        \begin{align}\label{eqLe4imp}
            \left|\pi_1\circ\Psi^{-1}(x,\tau)-x\right|&= \left| {\textstyle\sum}_{k=1}^2\varepsilon^k{v}_k\circ\Psi^{-1}(x,\tau)\right|.
        \end{align}
        From item \ref{lem:Lipschitz_conditions_2} in Lemma \ref{lem:Lipschitz_conditions} and the inequality \eqref{eq:lipschitz_near_identity}, we know that the maps $v_k\circ\Psi^{-1}$ are globally Lipschitz in $x$, uniformly in $\tau$, for $k\in\{1,2\}$, i.e., there exists constants $\hat{L}_{v,k} $ such that, for all $x_1,x_2\in\mathbb{R}^n$ and all $\tau\in\mathbb{R}_{\geq 0}$:
        \begin{align*}
            \left|v_k\circ\Psi^{-1}(x_1,\tau)-v_k\circ\Psi^{-1}(x_2,\tau)\right|&\leq \hat{L}_{v,k} |x_1-x_2|.
        \end{align*}
        Since the maps $v_k\circ\Psi^{-1}$ are also continuous and periodic in $\tau$, we may assume, without loss of generality, that, for all $(\tau,\varepsilon)\in\mathbb{R}_{\geq 0}\times[0,\bar{\varepsilon} ]$, we have $\left|v_k\circ\Psi^{-1}(0,\tau)\right|\leq \hat{L}_{v,k}$. Thus, adding and subtracting terms to \eqref{eqLe4imp}, and using the triangle inequality and the previous bounds, we obtain
        \begin{align*}
            \left|\pi_1\circ\Psi^{-1}(x,\tau)-x\right|&\leq  {\textstyle\sum}_{k=1}^2\varepsilon^k \hat{L}_{v,k} (|x|+1),
        \end{align*}
        for all $(x,\tau)\in\mathbb{R}^n\times\mathbb{R}_{\geq 0}$ and all $\varepsilon\in[0,\bar{\varepsilon} ]$. Using the reverse triangle inequality, we get
        \begin{align*}
            |x|-\varepsilon(|x|+1) {\textstyle\sum}_{k=1}^2 \hat{L}_{v,k}  \leq \left|\pi_1\circ\Psi^{-1}(x,\tau)\right|,
        \end{align*}
        for all $(x,\tau,\varepsilon)\in\mathbb{R}^n\times\mathbb{R}_{\geq 0}\times[0,\hat{\varepsilon}]$. Let $\hat{\varepsilon} =\min\{1,\bar{\varepsilon} ,1/((\hat{L}_{v,1} +\hat{L}_{v,2} ))\}$ and observe that, for all $ \varepsilon\in[0,\hat{\varepsilon} ]$, we have that $0\leq 1-\varepsilon {\textstyle\sum}_{k=1}^2 \hat{L}_{v,k} $
        for all $(x,\tau)\in\mathbb{R}^n\times\mathbb{R}_{\geq 0}$. Hence, we obtain that
        \begin{align*}
            \left(1-\varepsilon{\textstyle\sum}_{k=1}^2 \hat{L}_{v,k} \right)|x|-\varepsilon{\textstyle\sum}_{k=1}^2 \hat{L}_{v,k}  \leq \left|\pi_1\circ\Psi^{-1}(x,\tau)\right|,
        \end{align*}
        for all $(x,\tau,\varepsilon)\in\mathbb{R}^n\times\mathbb{R}_{\geq 0}\times[0,\hat{\varepsilon} ]$, which implies that
        \begin{align*}
            \left(1-\varepsilon{\textstyle\sum}_{k=1}^2 \hat{L}_{v,k} \right)\delta_3-\varepsilon{\textstyle\sum}_{k=1}^2 \hat{L}_{v,k}  \leq \left|\pi_1\circ\Psi^{-1}(x,\tau)\right|,
        \end{align*}
        for all $(x,\tau,\varepsilon)\in\mathcal{M}_3\times\mathbb{R}_{\geq 0}\times[0,\hat{\varepsilon} ]$. The result follows by defining $\varepsilon_0:=\min\{\hat{\varepsilon} ,(\delta_3-\delta_2)/((\delta_3+1)(\hat{L}_{{v},1} +\hat{L}_{{v},2} ))\}$. 
    \hfill$\blacksquare$ 
    \begin{lem}\thlabel{lem:avg_vector_field_structure}
        There exists a $\mathcal{C}^0$ map $g:\mathbb{R}^n\times\mathbb{R}_{\geq 0}\times[0,\varepsilon_0 ]\to\mathbb{R}^n$ such that, for all $(x,\tau)\in\mathcal{M}_3\times\mathbb{R}_{\geq 0}$ and all $\varepsilon\in[0,\varepsilon_0]$, the map $\Psi_*f_\varepsilon$, given by \eqref{eq:push_forward_map_definition}, satisfies
        \begin{align*}
            \Psi_*f_\varepsilon(x,\tau)&= \bar{f}(x) + \varepsilon\,g(x,\tau,\varepsilon),
        \end{align*}
        where $\bar{f}$ is given by \eqref{eq:averaged_vector_field}.
    \end{lem}
    \noindent\textbf{Proof:} By direct computation, we obtain:
        \begin{align*}
            \Psi_*f_\varepsilon(x,\tau)= &\left(\D x\Phi\circ\Psi^{-1}(x,\tau)\right)f_2\circ\Psi^{-1}(x,\tau)\\
            &~+\left(\D x\Phi\circ\Psi^{-1}(x,\tau)\right)f_1\circ\Psi^{-1}(x,\tau)\,\varepsilon^{-1}\\
            &~+\D \tau\Phi\circ\Psi^{-1}(x,\tau)\,\varepsilon^{-2}, 
        \end{align*}
        where
        \begin{align*}
            \D \tau\Phi\circ\Psi^{-1}(x,\tau)&= -  {\textstyle\sum}_{k=1}^2\varepsilon^k\D \tau{v}_k\circ\Psi^{-1}(x,\tau)\\
            \D x\Phi\circ\Psi^{-1}(x,\tau)&= \left(I -  {\textstyle\sum}_{k=1}^2\varepsilon^k\D x{v}_k\right)\circ\Psi^{-1}(x,\tau).
        \end{align*}
        Moreover, note that
        \begin{align*}
            &\D \tau{v}_1\circ\Psi^{-1}(x,\tau)= \hat{f}_1\circ\Psi^{-1}(x,\tau),
        \end{align*}
        and also that
        \begin{align*}
        \text{D}_{\tau}v_2(x,\tau)&=\underbrace{\hat{f}_2(x,\tau) + \text{D}_{x}\hat{f}_1(x,\tau) v_1(x,\tau) - \bar{f}(x)}_{=\D{\tau}w(x,\tau)}\\
        &\underbrace{- \text{D}_{x}\hat{f}_1(x,\tau) v_1(x,\tau) - \text{D}_{x}v_1(x,\tau) \hat{f}_1(x,\tau)}_{=-\D{\tau}\left(\D{x}v_1(x,\tau) v_1(x,\tau)\right)}\\
        &=\hat{f}_2(x,\tau)  - \text{D}_{x}v_1(x,\tau) \hat{f}_1(x,\tau) - \bar{f}(x)
        \end{align*}
        which implies that
        \begin{align*}
            \text{D}_{\tau}v_2\circ\Psi^{-1}(x,\tau)&= \hat{f}_2\circ\Psi^{-1}(x,\tau)- \bar{f}\circ\pi_1\circ\Psi^{-1}(x,\tau)\\
            &- \left(\text{D}_{x}v_1\circ\Psi^{-1}(x,\tau)\right)\hat{f}_1\circ\Psi^{-1}(x,\tau).
        \end{align*}
        Therefore, another direct computation shows that
        \begin{align}\label{eq:pushcan}
        \Psi_*f_\varepsilon(x,\tau)=&\left(f_1\circ\Psi^{-1}(x,\tau)-\hat{f}_1\circ\Psi^{-1}(x,\tau)\right)\varepsilon^{-1}\notag\\
            &~+\left(f_2\circ\Psi^{-1}(x,\tau)-\hat{f}_2\circ\Psi^{-1}(x,\tau)\right)\notag\\
            &~-\D x v_1\circ\Psi^{-1}(x,\tau) f_1\circ\Psi^{-1}(x,\tau)\notag\\
            &~+\D x v_1\circ\Psi^{-1}(x,\tau)\hat{f}_1\circ\Psi^{-1}(x,\tau)\notag\\
            &~+\bar{f}(x) + \hat{f}(x,\tau,\varepsilon),
        \end{align}
        where the map $\hat{f}$ is given by
        \begin{align*}
            \hat{f}(x,\tau,\varepsilon)&= \bar{f}\circ\pi_1\circ\Psi^{-1}(x,\tau)  - \bar{f}(x)\\
            &-\varepsilon   \D x {v}_1\circ\Psi^{-1}(x,\tau)\,f_2\circ\Psi^{-1}(x,\tau)\\
            &-\varepsilon   \D x {v}_2\circ\Psi^{-1}(x,\tau)\, f_1\circ\Psi^{-1}(x,\tau)\\
            &-\varepsilon^2 \D x {v}_2\circ\Psi^{-1}(x,\tau)\, f_2\circ\Psi^{-1}(x,\tau).
        \end{align*}
        Using Hadamard's Lemma \cite[Lemma 2.8]{nestruev2003smooth} and the fact that $\bar{f}$ is $\mathcal{C}^1$, we obtain:
        \begin{gather*}
            \bar{f}(x_1)  - \bar{f}(x_2)= \bar{F}(x_1,x_2)(x_1-x_2),
        \end{gather*}
        for all $x_1,x_2\in\mathbb{R}^n$, where $\bar{F}$ is given by
        \begin{equation*}
        \bar{F}(x_1,x_2):= {\textstyle\int}_0^1\D x \bar{f}(\lambda x_1+(1-\lambda) x_2)\,d\lambda.
        \end{equation*}
        Hence, using the fact that $\Psi^{-1}$ is a bijection: 
        \begin{align*}
            \bar{f}\circ\pi_1\circ\Psi^{-1}(x,\tau)  - \bar{f}(x)&= \tilde{F}(x,\tau)\,\left(\pi_1\circ\Psi^{-1}(x,\tau)-x\right),
        \end{align*}
        for all $x\in\mathbb{R}^n$, where $\tilde{F}$ is given by
        \begin{align}
            \tilde{F}(x,\tau):=\bar{F}\left(\pi_1\circ\Psi^{-1}(x,\tau),x\right).
        \end{align}
        However, since $\Psi$ is a diffeomorphism, we have
        \begin{align*}
            x=\pi_1\circ\Psi\circ\Psi^{-1}(x,\tau)=\pi_1\circ\Psi^{-1}(x,\tau) -  {\textstyle\sum}_{k=1}^2\varepsilon^k{v}_k\circ\Psi^{-1}(x,\tau),
        \end{align*}
        which implies that
        \begin{align*}
            \bar{f}\circ\pi_1\circ\Psi^{-1}(x,\tau)  - \bar{f}(x)&= \varepsilon\hphantom{^2} \tilde{F}\left(x,\tau\right){v}_1\circ\Psi^{-1}(x,\tau)\\
            &+ \varepsilon^2 \tilde{F}\left(x,\tau\right){v}_2\circ\Psi^{-1}(x,\tau),
        \end{align*}
        and that $\hat{f}$ can be written as
        \begin{align*}
            \hat{f}(x,\tau,\varepsilon)&= \varepsilon\,g(x,\tau,\varepsilon),
        \end{align*}
        where $g$ can be written in compact form as:
        \begin{align}\label{eqconstructiong}
            g&= \tilde{F}\,{v}_1\circ\Psi^{-1} -\D x{v}_1\circ\Psi^{-1}\,f_2\circ\Psi^{-1} - \D x{v}_2\circ\Psi^{-1}\,f_1\circ\Psi^{-1}\notag\\
            &+\varepsilon\,\left(\tilde{F}\,{v}_2\circ\Psi^{-1}-\D x{v}_2\circ\Psi^{-1}\,f_2\circ\Psi^{-1}\right).
        \end{align}
        Since $\Psi$ is a $\mathcal{C}^1$ diffeomorphism, and $g$ is a combination of $\mathcal{C}^0$ maps composed with $\Psi$, it follows that $g$ is $\mathcal{C}^0$ in all arguments. 
        
        Finally, by item \ref{lem:near_identity_transformation_properties_4} in Lemma \ref{lem:near_identity_transformation_properties}, for all $(x,\tau,\varepsilon)\in\mathcal{M}_3\times\mathbb{R}_{\geq 0}\times[0,\varepsilon_0 ]$, we have that
        $\Psi^{-1}(x,\tau)\in\mathcal{M}_2\times\mathbb{R}_{\geq 0}$. Also, by item \ref{lem:vector_fields_condition_time_varying_5} in Lemma \ref{lem:vector_fields_condition_time_varying}, for all $(x,\tau)\in\mathcal{M}_2\times\mathbb{R}_{\geq 0}$, we have that $\hat{f}_k(x,\tau)=f_k(x,\tau)$. Therefore,
        \begin{align*}
            \hat{f}_k\circ\Psi^{-1}(x,\tau)&=f_k\circ\Psi^{-1}(x,\tau),
        \end{align*}
        for all $(x,\tau,\varepsilon)\in\mathcal{M}_3\times\mathbb{R}_{\geq 0}\times[0,\varepsilon_0 ]$.
        Hence, in this set the first four terms in \eqref{eq:pushcan} cancel, and we obtain that pushforward map $\Psi_*f_\varepsilon$ satisfies
        \begin{equation*}
            \Psi_*f_\varepsilon({x},\tau)=\bar{f}({x}) + \varepsilon\,g({x},\tau,\varepsilon). 
        \end{equation*}
        for all $(x,\tau,\varepsilon)\in\mathcal{M}_3\times\mathbb{R}_{\geq 0}\times[0,\varepsilon_0 ]$. \hfill $\blacksquare$
    \begin{lem}\thlabel{lem:avg_vector_field_remainder_bounds}
        There exists a positive constant $L_g >0$ such that the map $g$, defined in \eqref{eqconstructiong}, satisfies
        \begin{align}
            |g(x,\tau,\varepsilon)|&\leq L_g (|{x}|+1),
        \end{align}
        for all $({x},\tau,\varepsilon)\in\mathbb{R}^n\times\mathbb{R}_{\geq 0}\times [0,\varepsilon_0 ]$.
    \end{lem}
    \noindent\textbf{Proof: }
        The map $g$ can be written in compact form as
        \begin{align}\label{eqexpansion}
            g(x,\tau,\varepsilon)&=  {\textstyle\sum}_{k=1}^5 G_i(x,\tau,\varepsilon)\,g_i(x,\tau,\varepsilon),
        \end{align}
        with the matrix-valued maps $G_i$ given by
        \begin{align*}
            G_1(x,\tau)&=\tilde{F}\left(x,\tau\right), & G_2(x,\tau)&= \D x{v}_1\circ \Psi^{-1}(x,\tau), \\
            G_3(x,\tau)&=\D x{v}_2\circ\Psi^{-1}(x,\tau), &
            G_4(x,\tau)&=\varepsilon\tilde{F}\left(x,\tau\right), \\
            G_5(x,\tau)&= \varepsilon \D x{v}_2\circ\Psi^{-1}(x,\tau),
        \end{align*}
        and the maps $g_i$ given by
        \begin{align*}
            g_1(x,\tau)&={v}_1\circ\Psi^{-1}(x,\tau), & g_2(x,\tau)&= -f_2\circ\Psi^{-1}(x,\tau), \\
            g_3(x,\tau)&=-f_1\circ\Psi^{-1}(x,\tau), & g_4(x,\tau)&={v}_2\circ\Psi^{-1}(x,\tau), \\
            g_5(x,\tau)&= -f_2\circ\Psi^{-1}(x,\tau).
        \end{align*}
        where the explicit (smooth) dependence on $\varepsilon$ is ommited to simplify notation. By Lemma \ref{lem:Lipschitz_conditions}, the maps $\bar{f}$ and $v_k$ are $\mathcal{C}^1$ and globally Lipschitz in $x$, uniformly in $\tau$. It follows that there exists constants $M_{g,i}>0$, $i\in\{1,2,\ldots,5\}$, such that, for all $(x,\tau,\varepsilon)\in\mathbb{R}^n\times\mathbb{R}_{\geq 0}\times[0,\varepsilon_0]$, we have that $|G_i(x,\tau,\varepsilon)|\leq M_{g,i}$ for all $i$. By Lemma \ref{lem:near_identity_transformation_properties}, the diffeomorphism $\Psi$ and its inverse $\Psi^{-1}$ are globally Lipschitz in $x$. In addition, from items \ref{asmp:vector_fields_condition_time_varying_1} in \thref{asmp:vector_fields_condition_time_varying} and Lemma \ref{lem:Lipschitz_conditions}, the maps $f_k$ are Lipschitz in $x$ for all $(x,\tau)\in\mathcal{M}_1\times\mathbb{R}_{\geq 0}$ and are $\mathcal{C}^0$ for all $(x,\tau)\in\mathbb{R}^n\times\mathbb{R}_{\geq 0}$. Since the composition of (globally Lipschitz) continuous functions is (globally Lipschitz) continuous, there exist constants $L_{g,i} >0$ such that, for all $x_1,x_2\in\mathcal{M}_3$ and for all $(\tau,\varepsilon)\in\mathbb{R}_{\geq 0}\times[0,\varepsilon_0 ]$, we have that $|g_i(x_1,\tau,\varepsilon)-g_i(x_2,\tau,\varepsilon)|\leq L_{g,i} |x_1-x_2|$. 
        
        Finally, since $g_i$ is continuous in all arguments, periodic in $\tau$, $[0,\varepsilon_0 ]$ is compact, and $\delta_3\mathbb{B}$ is compact, there exist constants $ M_{0,i}\in\mathbb{R}_{>0}$ such that $|g_i(x,\tau,\varepsilon)|\leq M_{0,i} $, for all $(x,\tau,\varepsilon)\in\delta_3\mathbb{B}\times\mathbb{R}_{\geq 0}\times[0,\varepsilon_0 ]$. Let $x_m\in\{x:\mathbb{R}^n:\,|x|=\delta_3\}$ be an arbitrary point, and note that $x_m\in\mathcal{M}_3\cap\delta_3\mathbb{B}$, and that, for all $(x,\tau,\varepsilon)\in\mathbb{R}^n\times\mathbb{R}_{\geq 0}\times[0,\varepsilon_0 ]$, we have that each term in \eqref{eqexpansion} can be written as:
        \begin{align*}
            G_i({x},\tau,\varepsilon) g_i({x},\tau,\varepsilon)&= G_i({x},\tau,\varepsilon)(g_i({x},\tau,\varepsilon)-g_i(x_m,\tau,\varepsilon))\\
            &+G_i({x},\tau,\varepsilon)g_i(x_m,\tau,\varepsilon).
        \end{align*}
        For all $x\in\mathbb{R}^n$, either $x\in\mathcal{M}_3$ or $x\in (\mathbb{R}^n\backslash\mathcal{M}_3)\subset\delta_3\mathbb{B}$. If $x\in\mathcal{M}_3$, then we have that
        \begin{align*}
            \left|G_i({x},\tau,\varepsilon) g_i({x},\tau,\varepsilon)\right|&\leq M_{g,i}L_{g,i}|x-x_m| + M_{g,i}M_{g,0}\\
            &\leq M_{g,i}L_{g,i}|x|+M_{g,i}(L_{g,i}\delta_3+M_{g,0}).
        \end{align*}
        Alternatively, if $x\in (\mathbb{R}^n\backslash\mathcal{M}_3)$, then we have that
        \begin{align*}
            \left|G_i({x},\tau,\varepsilon) g_i({x},\tau,\varepsilon)\right|\leq M_{g,i}M_{g,0}.
        \end{align*}
        Combining all of the above, we obtain that, for all $(x,\tau,\varepsilon)\in\mathbb{R}^n\times\mathbb{R}_{\geq 0}\times[0,\varepsilon_0 ]$, the map $g$ satisfies the inequality
        \begin{align*}
            |g(x,\tau,\varepsilon)|\leq L_{g} (|{x}| + 1), 
        \end{align*}
        where $L_g :=\max\left\{{\textstyle\sum}_{i=1}^5 M_{g,i} L_{g,i} ,{\textstyle\sum}_{i=1}^5 M_{g,i} (L_{g,i}\delta_3+M_{g,0}) \right\}$. \hfill $\blacksquare$ 
        
    \vspace{0.1cm}    
    All the claims of \thref{prop:near_identity_transformation} follow now directly by Lemmas \ref{lem:vector_fields_condition_time_varying}-\ref{lem:avg_vector_field_remainder_bounds}. \hfill$\blacksquare$
    \subsection{Proof of \thref{thm:GUES_implies_GUPES_timevarying}}
    Let Assumption \ref{asmp:vector_fields_condition_time_varying} generate $\delta_1$, let $\epsilon>0$, and let $\delta_2$ and $\delta_3$ satisfy \eqref{eq:deltas}. Let \thref{prop:near_identity_transformation} generate $\varepsilon_0 \in\mathbb{R}_{>0}$ such that its conclusions hold. Then, from \thref{asmp:vector_fields_condition_time_varying} and \thref{prop:near_identity_transformation}, and for all $\varepsilon\in(0,\varepsilon_0 )$, the map $\Psi_*f_\varepsilon$ is continuous. Hence, for all $(x_0,\tau_0,\varepsilon)\in\mathbb{R}^n\times\mathbb{R}_{\geq 0}\times(0,\varepsilon_0 )$, a solution to system \eqref{eq:perturbed_avged_system_timevaring} starting at the initial condition $(x_0,\tau_0)$ exists. 
    
    Let $V$ be given by \thref{asmp:Lyapunov_condition_timevarying}. Its time derivative along the trajectories of \eqref{eq:perturbed_avged_system_timevaring} satisfies:
    \begin{align*}
        \dot{V} &= \nabla V(x) ^\top\Psi_*f_\varepsilon(x,\tau,\varepsilon).
    \end{align*}
    Using item \ref{prop:near_identity_transformation_5} in \thref{prop:near_identity_transformation} and the bounds from \thref{asmp:Lyapunov_condition_timevarying}, we obtain that for all $(x,\tau,\varepsilon)\in\mathcal{M}_3\times\mathbb{R}_{\geq 0}\times(0,\varepsilon_0 )$:
    \begin{align*}
        \dot{V}&\leq \nabla V(x)^\top\bar{f}({x}) + \varepsilon \nabla V(x)^\top g({x},\tau,\varepsilon)\\
        &\leq -c_1 \phi({x})^2 + \varepsilon c_2 \phi({x})|g(x,\tau,\varepsilon)|.
    \end{align*}
    We consider two possible cases:
    \begin{enumerate}[C1),topsep=0pt,itemsep=-1ex,partopsep=1ex,parsep=1ex]
        \item Item \ref{asmp:perturbation_norm_upper_bound}-\ref{asmp:perturbation_norm_upper_bound_1} in \thref{asmp:Lyapunov_condition_timevarying} holds, and in this case we obtain:
        \begin{align*}
            \dot{V}\leq &- \frac{c_1}{2}\phi({x})^2 - \left(\frac{c_1}{4}-\varepsilon c_2 \bar{L}_g \right) \phi({x})^2\\
            &-\phi(x)\left(\frac{c_1}{4}\phi(x) -\varepsilon c_2 \bar{L}_g \right),
        \end{align*}
        for all $(x,\tau,\varepsilon)\in\mathcal{M}_3\times\mathbb{R}_{\geq 0}\times(0,\varepsilon_0 )$. Since the function $\phi$ is positive definite, there exists a class $\mathcal{K}$ function $\alpha_4$ such that, for all $x\in\mathbb{R}^n$, $\alpha_4(|x|)\leq \phi(x)$ \cite[Lemma 4.3]{khalil2002nonlinear} and $\alpha_4^{-1}(s)$ is defined for $s$ sufficiently small. Therefore, there exists $\varepsilon^*_1 \in\left(0,\min\left\{\varepsilon_0,\eta^*\right\}\right)$, where $\eta^*=\frac{c_{1}}{4c_2\bar{L}_g}$, such that, for all $\varepsilon\in(0,\varepsilon^*_1 )$: 
        \begin{align*}
            \dot{V}&\leq - \frac{c_1}{2}\phi({x})^2, & \forall |x|\geq \delta_3 + \alpha_4^{-1}\left(\frac{\varepsilon}{\eta^*}\right).
        \end{align*}
        \item Item \ref{asmp:perturbation_norm_upper_bound}-\ref{asmp:perturbation_norm_upper_bound_2} in \thref{asmp:Lyapunov_condition_timevarying} holds, and in this case:
        \begin{align*}
            \dot{V}
            \leq -\frac{c_1}{2}\phi({x})^2 &- \phi(x)\left(\frac{c_1}{4}\alpha_3(|{x}|)-\varepsilon c_2 L_g \right)|{x}|\\
            &-\phi({x})\left(\frac{c_1}{4}\alpha_3(|{x}|)|{x}|-\varepsilon c_2  L_g \right),
        \end{align*}
        for all $(x,\tau,\varepsilon)\in\mathcal{M}_3\times\mathbb{R}_{\geq 0}\times(0,\varepsilon_0 )$, where $\alpha_3\in\mathcal{K}$, and where we used  \thref{prop:near_identity_transformation}-(e). Thus, there exists $\varepsilon_2^* \in(0,\varepsilon_0 )$ such that, for all $\varepsilon\in(0,\varepsilon_2^* )$:
        \begin{align*}
            \dot{V}\leq -\frac{c_1}{2}\phi(x)^2,~~\forall |{x}|\geq \delta_3+\alpha_3^{-1}\left(\frac{\varepsilon}{\eta^*}\right)+\alpha_5^{-1}\left(\frac{\varepsilon}{\eta^*}\right)
        \end{align*}
      where $\alpha_5(|x|):=\alpha_3(|x|)|x|$, and where, without loss of generality, we assumed $\bar{L}_g\geq L_g$.
    \end{enumerate}
    \vspace{0.1cm}
     Let $\varepsilon^* :=\min\{1,\varepsilon_1^* ,\varepsilon_2^* \}$,  $\rho(\varepsilon,\delta_1):=\sum_{i=1}^3\alpha_{i+2}^{-1}(\varepsilon/\eta^*(\delta_1))$. Then combining the above two cases we obtain that $\dot{V}\leq -c\phi(x)^2$, for all $|x|\geq 2\delta_3+\rho(\varepsilon,\delta_1)$ and all $\varepsilon\in(0,\varepsilon^*)$, with $c:=c_1/2$. Then, following similar steps as in \cite[Appendix C.9]{khalil2002nonlinear} and the proof of \cite[Appendix C.]{cai2009characterizations}, there exist functions  $\beta\in\mathcal{KL}$ and $\kappa\in\mathcal{K}_{\infty}$ such that, for all $(x_0,\tau_0,\varepsilon)\in\mathbb{R}^n\times\mathbb{R}_{\geq0}\times(0,\varepsilon^* )$, any solution to system \eqref{eq:perturbed_avged_system_timevaring} starting at $(x_0,\tau_0)$, satisfies
    \begin{align}\label{KLbound111}
        |{x}(t)|&\leq \beta(|{x}(t)|,t) + \Delta_{\delta,\varepsilon}, ~~~~\forall t\geq0,
    \end{align}
    where $\Delta_{\delta,\varepsilon}:=\kappa(2\delta_3+\rho(\varepsilon,\delta_1))$. 
\hfill$\blacksquare$
 \subsection{Proof of \thref{cor:GUES_implies_GUPES_timevarying_1}}
    By \thref{thm:GUES_implies_GUPES_timevarying}, there exists an $\varepsilon^* \in(0,\varepsilon_0 )$ and $\beta\in\mathcal{KL}$ such that, for all $(\bar{x}_0,\tau_0,\varepsilon)\in\mathbb{R}^n\times\mathbb{R}_{\geq0}\times(0,\varepsilon^* )$, any solution $(\bar{x},\tau)$ of system \eqref{eq:perturbed_avged_system_timevaring} starting from  $(\bar{x}_0,\tau_0)$ satisfies $|\bar{x}(t)|\leq\beta(|\bar{x}_0|,t) + \Delta(\varepsilon,\delta)$, for all $t\geq0$, with  $\Delta(\varepsilon,\delta):=\kappa(\rho(\delta,\varepsilon))$. For any $(x_0,\tau_0)$, let $(x,\tau)$ be a solution of system \eqref{eq:orig_system_timevaring} starting from the initial condition $(x_0,\tau_0)$. Since $\Psi$ is a $\mathcal{C}^1$ diffeomorphism, and system \eqref{eq:perturbed_avged_system_timevaring} is the pushforward of system \eqref{eq:orig_system_timevaring} under $\Psi$, it follows that $(x(t),\tau(t))=\Psi^{-1}(\bar{x}(t),\tau(t))$, for all $t>0$, where $(\bar{x},\tau)$ is a solution of system \eqref{eq:perturbed_avged_system_timevaring} with initial condition $\Psi(x_0,\tau_0)$. Therefore, for all $t\geq0$ we have:
    \begin{align*}
        |x(t)|&=\left|\pi_1\circ\Psi^{-1}(\bar{x}(t),\tau(t))\right|\\
        &\leq\left|\pi_1\circ\Psi^{-1}(\bar{x}(t),\tau(t))-\pi_1\circ\Psi^{-1}(0,\tau(t))\right|\\
        &+\left|\pi_1\circ\Psi^{-1}(0,\tau(t))\right|.
    \end{align*}
    From item 
    \ref{prop:near_identity_transformation_3} in \thref{prop:near_identity_transformation}, we obtain that
    \begin{align*}
        \left|\pi_1\circ\Psi^{-1}(\bar{x}(t),\tau(t))-\pi_1\circ\Psi^{-1}(0,\tau(t))\right|\leq (1+L_{\Psi} \varepsilon)|\bar{x}(t)|,
    \end{align*}
    and also $|\pi_1\circ\Psi^{-1}(0,\tau(t))|\leq L_\Psi  \varepsilon$. Therefore, it follows that
    \begin{align*}
        |x(t)|&\leq (1+L_\Psi \varepsilon) \beta(|\pi_1\circ \Psi(x_0,\tau_0)|,t)\\
        &~~~~+ (1+L_\Psi \varepsilon) \Delta_{\delta,\varepsilon} + L_\Psi \varepsilon.
    \end{align*}
    Similarly, we have
    \begin{align*}
        |\pi_1\circ \Psi(x_0,\tau_0)|&\leq |\pi_1\circ \Psi(x_0,\tau_0)-\pi_1\circ \Psi(0,\tau_0)|+|\pi_1\circ \Psi(0,\tau_0)|\\
        &\leq (1+L_\Psi \varepsilon)|x_0| + L_\Psi \varepsilon.
    \end{align*}
    Since $\beta(\cdot,t)\in\mathcal{K}_\infty$, it is strictly increasing and satisfies
    \begin{align*}
        \beta(|\pi_1\circ \Psi(x_0,\tau_0)|,t)&\leq \beta((1+L_\Psi \varepsilon)|x_0| + L_\Psi \varepsilon,t).
    \end{align*}
    for all $(x_0,\tau_0,t,\varepsilon)\in\mathbb{R}^n\times\mathbb{R}_{\geq 0}\times\mathbb{R}_{\geq0}\times[0,\varepsilon^* ]$. We then have two possible cases:
    \vspace{0.1cm}
    \begin{enumerate}[C1),topsep=0pt,itemsep=-1ex,partopsep=1ex,parsep=1ex]
        \item If $|x_0|\leq L_\Psi\varepsilon$, then
        \begin{align*}
            \beta(|\pi_1\circ \Psi(x_0,\tau_0)|,t)&\leq \beta(2L_\Psi \varepsilon+L_\Psi^2\varepsilon^2 ,t).
        \end{align*}
        %
        %
        \item If $|x_0| >L_\Psi\varepsilon$, then
        \begin{align*}
            \beta(|\pi_1\circ \Psi(x_0,\tau_0)|,t)\leq \beta((2+L_\Psi\varepsilon )|x_0|,t).
        \end{align*}
    \end{enumerate}
    Therefore, for all $(x_0,\tau_0,t,\varepsilon)\in\mathbb{R}^n\times\mathbb{R}_{\geq 0}\times\mathbb{R}_{\geq 0}\times[0,\varepsilon^*]$:
    \begin{align*}
        \beta(|\pi_1\circ \Psi(x_0,\tau_0)|,t)&\leq \beta((2+L_\Psi\varepsilon )|x_0|,t)+\beta(2L_\Psi \varepsilon+L_\Psi^2\varepsilon^2,t)\\
        &\leq \beta((2+L_\Psi\varepsilon)|x_0|,t)+\beta(2L_\Psi \varepsilon+L_\Psi^2\varepsilon^2,0).
    \end{align*}
    However, from Claim \ref{refclaim1}, we have that $L_\Psi\varepsilon\leq 1$, for all $\varepsilon\in[0,\varepsilon^*]$. Therefore, we have that
    \begin{align*}
        \beta(|\pi_1\circ \Psi(x_0,\tau_0)|,t)&\leq \beta(3|x_0|,t)+\beta(3L_\Psi\varepsilon,0).
    \end{align*}
    The result of the corollary follows by defining
    \begin{equation}\label{eq:Delta_tilde}
    \begin{aligned}
            \tilde{\Delta}(\delta,\varepsilon)&:=(1+L_\Psi \varepsilon)\Delta_{\delta,\varepsilon}+L_\Psi \varepsilon + (1+L_\Psi \varepsilon)\beta(3L_\Psi\varepsilon,0)
    \end{aligned}
    \end{equation}
    and 
    \begin{equation}\label{KLtilde}
    \tilde{\beta}(r,s):=2\beta(3 r,s),
    \end{equation}
    for all $r,s\geq0$. \hfill$\blacksquare$
    \subsection{Proof of \thref{cor:GUES_implies_GUPES_timevarying_2}}
    {It follows from the assumptions of the corollary that $\kappa\in\mathcal{K}$ and $\beta\in\mathcal{KL}$ generated by the proof of Theorem 1 are independent of the choice of $\delta_1>0$, $\epsilon>0$, and $\delta$ satisfying \eqref{eq:deltas}.} Let $\nu>0$ be given. Hence, there exists $r>0$ sufficiently small such that $\kappa(r)<\frac{v}{4}$. Let $\delta_1<r/16$, and choose $\epsilon=\sqrt{3}-1$, $\delta_2=(1+\epsilon)\delta_1$, $\delta_3=\frac{4}{3}(1+\epsilon)\delta_2=\frac{4}{3}(1+\epsilon)^2\delta_1=4\delta_1$, which satisfy \eqref{eq:deltas}. It follows that $2\delta_3=8\delta_1<\frac{r}{2}$. Let $\varepsilon_a>0$ be such that $\sum_{i=1}^3\alpha_{i+2}^{-1}(\varepsilon/(\eta^*(\delta_1)))<\frac{r}{2}$ for all $\varepsilon\in(0,\varepsilon_a)$. Such $\varepsilon_a$ always exist because $\alpha_{i+2}^{-1}\in\mathcal{K}$. It follows that $\Delta_{\delta,\varepsilon}=\kappa(2\delta_3+\sum_{i=1}^3\alpha_{i+2}^{-1}(\varepsilon/(\eta^*(\delta_1))))<\frac{v}{4}$. Let $\varepsilon_b>0$ be such that $L_\Psi \varepsilon \leq \min\{\nu/3,1/3\}$ and $\beta(3L_\Psi \varepsilon,0)\leq \nu/4$ for all $\varepsilon\in(0,\varepsilon_b)$, where $\beta$ comes from \eqref{KLbound111}. Such $\varepsilon_b$ always exist because $\beta(\cdot,s)\in\mathcal{K}$. Let $\varepsilon^*>0$ be generated by Corollary 1, and define $\varepsilon^{**}=\min\{\varepsilon_a,\varepsilon_b,\varepsilon^*\}$. Then, every solution of system \eqref{eq:orig_system_timevaring} starting at $(x_0,\tau_0)$ satisfies the bound \eqref{UUbound} with $\mathcal{K}\mathcal{L}$ function given by \eqref{KLtilde} and ultimate bound $\tilde{\Delta}_{\delta,\varepsilon}$ given by \eqref{eq:Delta_tilde}. However, by the choice of $\delta$ and $\varepsilon$, we have that $\tilde{\Delta}_{\delta,\varepsilon}\leq \nu$, which establishes the desired bound.\hfill  $\blacksquare$

    
    \subsection{Proof of \thref{cor:GUES_implies_GUPES_timevarying_3}}
    Since \thref{asmp:vector_fields_condition_time_varying} is satisfied for $\delta_1=0$, we may pick $\delta_1=\delta_2=0$, and $\delta_3\in(0,\infty)$ arbitrarily small. Following similar steps to the proof of \thref{cor:GUES_implies_GUPES_timevarying_2} yields the result. \hfill $\blacksquare$
    \subsection{Proof of \thref{thm:es_unbounded_gradient}}
    \label{sectionprooftheorem2}
        We first consider the case when $J^\star\in\mathbb{R}$ is arbitrary, and we verify that the maps defining system \eqref{eq:closed_loop_es_orig_sys} satisfy \thref{asmp:vector_fields_condition_time_varying}. Clearly, the right hand side in \eqref{eq:closed_loop_es_orig_sys} is $\mathcal{C}^0$ and satisfies item \ref{asmp:vector_fields_condition_time_varying_2} in \thref{asmp:vector_fields_condition_time_varying}. Let $\bar{J}\in\mathbb{R}_{>0}$, and let $\delta_1\in[0,\infty)$ be such that $J(x)\geq \bar{J}$, for all $|x|\geq \delta_1$. Such  $\delta_1$ always exists because $J$ is radially unbounded. It follows that the feedback law
        \begin{align*}
            u_{i,1}(J(x),\tau)&=\sqrt{2\omega_i J(x)}\cos(\log(J(x))+\omega_i\tau)\\
            u_{i,2}(J(x),\tau)&=\sqrt{2\omega_i J(x)}\sin(\log(J(x))+\omega_i\tau),
        \end{align*}
        is $\mathcal{C}^2$ for all $|x|\geq \delta_1$. Therefore, system \eqref{eq:closed_loop_es_orig_sys} satisfies item \ref{asmp:vector_fields_condition_time_varying_3} in \thref{asmp:vector_fields_condition_time_varying}. Next, via direct differentiation, we obtain that
        \begin{align*}
            \left|\D{x}\left(u_{i,j}^{\varepsilon,1}(J(x),\tau)\right)\right|&\leq \sqrt{2\omega_i}J(x)^{-\frac{1}{2}}|\nabla J(x)|.
        \end{align*}
        From \thref{asmp:cost_function_conditions_es_1}-(c), we have that $\nabla J$ is $L_J$-globally Lipschitz, which implies that $|\nabla J(x)^2|\leq 2L_J(J(x)-J^*)$, for all $x\in\mathbb{R}^n$ \cite[Lemma 1, pp.23]{PolyakBook}. It follows that 
        \begin{align}\label{boundonjacobiantheorem2}
            \left|\D{x}\left(u_{i,1}^{\varepsilon,1}(J(x),\tau)\right)\right|&\leq \sqrt{2\omega_i\kappa}(1-J^\star J(x)^{-1})^{\frac{1}{2}}\leq m,
        \end{align}
        for all $|x|\geq \delta_1$ and all $\tau\in\mathbb{R}_{\geq 0}$, where $m:=\sqrt{2\omega_i\kappa}(1+|J^\star|\bar{J}^{-1})^{\frac{1}{2}}$ and $\kappa=2L_J$.
        Therefore, system \eqref{eq:closed_loop_es_orig_sys} satisfies item \ref{asmp:vector_fields_condition_time_varying_1} in \thref{asmp:vector_fields_condition_time_varying}. A similar computation shows that system \eqref{eq:closed_loop_es_orig_sys} satisfies item \ref{asmp:vector_fields_condition_time_varying_4} in \thref{asmp:vector_fields_condition_time_varying}. 
        
        Next, let $\epsilon\in(0,\infty)$,  $\delta_2\in[(1+\epsilon)\delta_1,\infty)$, $\delta_3\in((1+\epsilon)\delta_2,\infty)$, and let $\mathcal{M}_j$, for $j\in\{1,2,3\}$, be the corresponding nested subsets defined in \eqref{eq:concentric_subsets}. Using the formula \eqref{eq:averaged_vector_field}, the nominal average system \eqref{eq:nominal_avged_sys} corresponding to system \eqref{eq:closed_loop_es_orig_sys} on $\mathcal{M}_3$, is given by
        \begin{align}\label{eq:avged_sys_es}
            \dot{\tilde{x}}&=\bar{f}(\tilde{x})=b_0(\tilde{x}+x^\star) + \textstyle\sum_{\substack{i=1}}^{r}\sum_{j=1}^2 b_{i,j}b_{i,j}^\top \nabla J(\tilde{x}+x^\star).
        \end{align}
        Consider the Lyapunov function candidate $V$ and the positive definite function $\phi$ defined by
        \begin{align}\label{lyapunov001}
        V(\tilde{x})&=J(\tilde{x}+x^\star)-J(x^\star), & \phi(\tilde{x})&=|\nabla J(\tilde{x}+x^\star)|,
        \end{align}
        which satisfy the inequalities
        \begin{gather}\label{lyapunov002}
            \alpha_1(|\tilde{x}|)\leq V(\tilde{x})\leq \alpha_2(|\tilde{x}|),\\
            \label{lyapunov003}|\nabla V(\tilde{x})|\leq \phi(\tilde{x}),
        \end{gather}
        for all $x\in\mathbb{R}^n$, and satisfy the inequality
        \begin{gather}
            \label{lyapunov004}\nabla V(\tilde{x})^\top\bar{f}(\tilde{x})\leq (\kappa_3-\gamma)|\nabla J(\tilde{x}+x^\star)|^2 < 0
        \end{gather}
        for all $x\in\mathcal{M}_3$, where the functions $\alpha_1$ and $\alpha_2$ are $\mathcal{K}_\infty$ functions, whose existence is guaranteed by the radial unboundedness of $V$ \cite[Lemma 4.3]{khalil2002nonlinear}. Since, by assumption $\gamma>\kappa_3$, system \eqref{eq:closed_loop_es_orig_sys} satisfies items \ref{asmp:Lyapunov_condition_timevarying_1}-\ref{asmp:Lyapunov_condition_timevarying_2} in \thref{asmp:Lyapunov_condition_timevarying}. Moreover, from \thref{asmp:cost_function_conditions_es_1}, we have that $\alpha_J(|\tilde{x}|)|\tilde{x}|\leq |\nabla J(\tilde{x}+x^\star)|$, where $\alpha_J$ is a class $\mathcal{K}$ function. Hence, system \eqref{eq:closed_loop_es_orig_sys} satisfies item \ref{asmp:perturbation_norm_upper_bound}-\ref{asmp:perturbation_norm_upper_bound_2} in \thref{asmp:Lyapunov_condition_timevarying} with $\phi(\tilde{x}):=|\nabla J(\tilde{x}+x^\star)|$. By \thref{thm:GUES_implies_GUPES_timevarying} we conclude that system \eqref{eq:closed_loop_es_orig_sys} is UGUB. 

        Next, we consider the case when $J^\star\in\mathbb{R}_{>0}$. In this case, there exists $\bar{J}\in\mathbb{R}$ such that $0<\bar{J}<J^\star$. Since $J^\star$ is the minimum value of the cost, it follows that $J(x)>\bar{J}$ for all $x\in\mathbb{R}^n$, which implies that the previous computations hold with $\delta_1=0$. In addition, in this case \thref{asmp:vector_fields_condition_time_varying} holds with $\delta_1=0$, {and the Lyapunov function candidate $V$ and the positive definite function $\phi$ in \eqref{lyapunov001} still satisfy the inequalities \eqref{lyapunov002}-\eqref{lyapunov004} for any choice of $\epsilon\in(0,\infty)$ and $\delta$ satisfying \eqref{eq:deltas}, with $\mathcal{M}_j$, for $j\in\{1,2,3\}$, being the corresponding nested subsets defined in \eqref{eq:concentric_subsets}}. Therefore, by invoking \thref{cor:GUES_implies_GUPES_timevarying_3}, we conclude that system \eqref{eq:closed_loop_es_orig_sys} is UGpAS. Finally, if $J^\star = 0$, we can take $\bar{J}=0$ and in this case, \thref{asmp:vector_fields_condition_time_varying} will be satisfied for all $\delta_1\in(0,\infty)$ using $\tilde{m}=\sqrt{2\omega_i\kappa}$ in \eqref{boundonjacobiantheorem2}. Therefore, by \thref{cor:GUES_implies_GUPES_timevarying_2}, we conclude that the closed-loop system is UGpAS. \hfill $\blacksquare$

    \subsection{Proof of \thref{thm:es_bounded_gradient}}
         It is easy to see that the right hand side in \eqref{eq:closed_loop_es_orig_sys} is $\mathcal{C}^0$ and satisfies item \ref{asmp:vector_fields_condition_time_varying_2} in \thref{asmp:vector_fields_condition_time_varying}. In addition, the maps $u_{i,j}$ are $\mathcal{C}^\infty$, which implies that the maps $u_{i,j}(J(\cdot),\tau)$ are $\mathcal{C}^2$ for all $x\in\mathbb{R}^n$. Therefore, system \eqref{eq:closed_loop_es_orig_sys} satisfies item \ref{asmp:vector_fields_condition_time_varying_3} in \thref{asmp:vector_fields_condition_time_varying}. Also, note that
        \begin{align}\label{eq:bound_on derivative_es_2}
            \left|\D{x}\left(u_{i,j}(J(x),\tau)\right)\right|&\leq \sqrt{2\omega_i}|\nabla J(x)|\leq \sqrt{2\omega_i} M_J,
        \end{align}
        for all $x\in\mathbb{R}^n$, where we used item b) in \thref{asmp:cost_function_conditions_es_1}. Therefore, system \eqref{eq:closed_loop_es_orig_sys} satisfies item \ref{asmp:vector_fields_condition_time_varying_1} in \thref{asmp:vector_fields_condition_time_varying}. A similar computation shows that system \eqref{eq:closed_loop_es_orig_sys} satisfies item \ref{asmp:vector_fields_condition_time_varying_4} in \thref{asmp:vector_fields_condition_time_varying}. Therefore, system \eqref{eq:closed_loop_es_orig_sys} satisfies \thref{asmp:vector_fields_condition_time_varying}. Next, let $\delta_1=\delta_2=0$, and fix a choice of $\delta_3\in(0,\infty)$, and let $\mathcal{M}_j$ for $j\in\{1,2,3\}$ be the corresponding nested subsets as defined in \eqref{eq:concentric_subsets}. Using the formula \eqref{eq:averaged_vector_field}, the nominal average system \eqref{eq:nominal_avged_sys} corresponding to system \eqref{eq:closed_loop_es_orig_sys} on $\mathbb{R}^n$, is given by
        \begin{align}
            \dot{\tilde{x}}&=\bar{f}(\tilde{x})=b_0(\tilde{x}+x^\star) + \textstyle\sum_{\substack{i=1}}^{r}\sum_{j=1}^2 b_{i,j}b_{i,j}^\top \nabla J(\tilde{x}+x^\star).        
        \end{align}
        Since the Lyapunov function candidate $V$ and the function $\phi$ \eqref{lyapunov001} now satisfy the inequalities \eqref{lyapunov002}-\eqref{lyapunov004} for all $x\in\mathcal{M}_3$, {for any choice of $\delta_3\in(0,\infty)$},
        it follows that system \eqref{eq:closed_loop_es_orig_sys} satisfies items \ref{asmp:Lyapunov_condition_timevarying_1}-\ref{asmp:Lyapunov_condition_timevarying_2} in \thref{asmp:Lyapunov_condition_timevarying}. Moreover, from item b) in \thref{asmp:cost_function_conditions_es_1}, we have that $\alpha_J(|\tilde{x}|)\leq |\nabla J(\tilde{x}+x^\star)|$, where $\alpha_J$ is a class $\mathcal{K}$ function. Finally, using item b) in \thref{asmp:cost_function_conditions_es_1} and item d) in \thref{asmp:drift_term_es}, we obtain that $|b_0(x)|\leq |\nabla J(x)|\leq L_J$, for all $x\in\mathbb{R}^n$. Then, we have the following Claim.
        \begin{claim}\thlabel{claim:boundedness_of_g}
            The remainder map $g$ from item \ref{prop:near_identity_transformation_5} in \thref{prop:near_identity_transformation} is uniformly bounded, for all $(x,\tau,\varepsilon)\in\mathbb{R}^n\times\mathbb{R}_{\geq 0}\times [0,\varepsilon_0]$.
        \end{claim}
        The proof of \thref{claim:boundedness_of_g} can be found in the Appendix. From \thref{claim:boundedness_of_g}, it follows that system \eqref{eq:closed_loop_es_orig_sys} satisfies item \ref{asmp:perturbation_norm_upper_bound}-\ref{asmp:perturbation_norm_upper_bound_1} in \thref{asmp:Lyapunov_condition_timevarying}. Therefore, by \thref{cor:GUES_implies_GUPES_timevarying_3} we  conclude that system \eqref{eq:closed_loop_es_orig_sys} is UGpAS. \hfill $\blacksquare$
\section{Conclusion and Future Work}\label{sec:conclusions}
In this manuscript, we introduced a (second-order) averaging method that allows us to study the stability properties of a class of oscillatory systems with periodic flows based on the stability properties of their corresponding averaged systems. In contrast to existing results in the literature, the method is suitable for the study of uniform global (practical) stability properties. Such properties are studied under suitable assumptions, which, naturally, are stronger compared to others that only enable local or semi-global practical results. By leveraging the proposed method, we showed that a class of extremum seeking algorithms is able to achieve uniform global practical asymptotic stability for a broad range of cost functions, which include quadratic (with positive definite Hessian), strongly convex, and certain invex functions. Future research will extend these results via singular perturbation theory to study dynamic plants in the loop, as well as systems with hybrid dynamics. 

\bibliographystyle{elsarticle-num}
\bibliography{global_feedback}

\appendix
\section{Auxiliary Lemmas}
\begin{lem}\thlabel{lem:strong_convexity_implies_assumption}
    Let $J:\mathbb{R}^n\rightarrow\mathbb{R}$ be a $\mu$-strongly convex $\mathcal{C}^1$ function with $L$ globally Lipschitz gradient. Then, item (d) in Assumption \ref{asmp:cost_function_conditions_es_1} is satisfied. 
\end{lem}
\textbf{Proof:} The upper bound follows directly by \cite[Thm. 2.1.5]{NesterovsBook}. To obtain the lower bound, not that by $\mu$-strong convexity:
\begin{align*}
    \left(\nabla J(x_1)-\nabla J(x_2)\right)^\top(x_1-x_2)\geq \mu |x_1-x_2|^2,
\end{align*}
for all $x_1,x_2\in\mathbb{R}^n$. Using Cauchy-Schwartz inequality, it is easy to see that the following holds  $\left|\nabla J(x_1)-\nabla J(x_2)\right|\geq \mu |x_1-x_2|$ for all $x_1,x_2\in\mathbb{R}^n$. It follows that
%
%
%
\begin{align*}
    \left|\nabla J(x_1)-\nabla J(x_2)\right|^2\geq \mu^2 |x_1-x_2|^2.
\end{align*}
for all $x_1,x_2\in\mathbb{R}^n$. Let $\alpha_J:[0,\infty)\rightarrow [0,\mu)$ be given by $\alpha_J(s):= \mu \,\tanh(s)$, which is strictly increasing and satisfies $\alpha_J(0)=0$. Therefore, $\alpha_J\in\mathcal{K}$ and, by definition, $\alpha_J(s)<\mu$, for all $s\geq0$. It follows that
\begin{align*}
    \left|\nabla J(x_1)-\nabla J(x_2)\right|^2\geq \alpha_J(|x_1-x_2|)^2|x_1-x_2|^2.
\end{align*}
for all $x_1,x_2\in\mathbb{R}^n$.
\begin{lem}\thlabel{lem:exmp_2_cost_function}
    Let $J:\mathbb{R}^n\rightarrow\mathbb{R}$ be the function defined in \thref{exmp:es_1_illustration}. Then, $J$ satisfies \thref{asmp:vector_fields_conditions_es}.
\end{lem}
\noindent\textbf{Proof:} The cost function $J$ can be written as $J=h\circ H$ where $h(s)=s+3\sin(\sqrt{s})^2$ and $H(x):=|x-x^\star|^2$ are $\mathcal{C}^\infty$ everywhere on their domain. Moreover, $H(x)\geq 0$ for all $x\in\mathbb{R}^n$. Therefore, the function $J=h\circ H$ is $\mathcal{C}^\infty$. The derivative of $J$ satisfies
\begin{align*}
    \nabla J(x) = \D{} h(H(x))\nabla H(x) = 2\D{} h(H(x))(x-x^\star),
\end{align*}
where
\begin{align*}
    \D{} h(H(x)) = \frac{1}{2}\left(2+\frac{3\sin\left(2\sqrt{H(x)}\right)}{\sqrt{H(x)}}\right)\in\mathbb{R}.
\end{align*}
It follows that
\begin{align*}
    |\nabla J(x)|^2 = 4|\D{} h(H(x))(x-x^\star)|^2=4Dh(H(x))^2(x-x^{\star})^2,
\end{align*}
and it can be verified that $\frac{1}{4}<\D{} h(H(x))<4$, for all $x\in\mathbb{R}^n$. Therefore, there exists $\mu\in\mathbb{R}_{>0}$ such that, for all $x\in\mathbb{R}^n$, we have that
\begin{align*}
    |\nabla J(x)|^2\geq 4\mu^2|x-x^\star|^2\geq \alpha_J(|x-x^\star|)^2|x-x^\star|^2,
\end{align*}
where $\alpha_J(s):=2\mu\tanh(s)$. 
Similarly, the second derivative of $J$ satisfies
\begin{align*}
    \nabla^2J(x)=\D{}^2h(H(x))\nabla H(x)\nabla H(x)^\top + \D{}h(H(x))\nabla^2H(x),
\end{align*}
where 
\begin{align*}
    \D{}^2h(H(x))= \frac{3 \cos \left(2 \sqrt{H(x)}\right)}{2 H(x)}-\frac{3 \sin \left(2 \sqrt{H(x)}\right)}{4 H(x)^{3/2}},
\end{align*}
It follows that
\begin{align*}
    |\nabla^2J(x)|&\leq |\D{}^2h(H(x))||\nabla H(x)|^2+|\D{}h(H(x))||\nabla^2 H(x)|,
\end{align*}
%
where $|\D{}^2h(H(x))|\leq \frac{3}{H(x)}$. Hence, the Hessian satisfies the inequality
\begin{align*}
    |\nabla^2J(x)|&\leq \frac{3|\nabla H(x)|^2}{H(x)}+8\leq 20.
\end{align*}
\section{Proofs of Auxiliary Claims}
\subsection{Proof of Claim \ref{refclaim1}}
\label{profclaim1}
\noindent\textbf{Proof: } The matrix $\D{x}\Phi(x,\tau)$ is a square matrix, and therefore its singular value decomposition is given by
        \begin{align*}
            \D{x}\Phi(x,\tau) = V(x,\tau)\Sigma(x,\tau) U(x,\tau)^\top,
        \end{align*}
        where the matrices $V(x,\tau)$ and $U(x,\tau)$ are orthonormal matrices and $\Sigma(x,\tau)$ is a square diagonal matrix with the singular values of $\D{x}\Phi(x,\tau)$ on the diagonal. Since for all $\varepsilon\in[0,\bar{\varepsilon}_2 ]$, for all $(x,\tau)\in\mathbb{R}^n\times\mathbb{R}_{\geq 0}$, the eigenvalues of the Jacobian matrix $\D{x}\Phi(x,\tau)$ are contained in the compact interval $[1-2\tilde{L}_{\Psi} \varepsilon,1+2\tilde{L}_{\Psi} \varepsilon]\subset [1/2,3/2]$, it follows that the singular values of $\D{x}\Phi(x,\tau)$ coincide with its eigenvalues and therefore are also contained in the compact interval $[1-2\tilde{L}_{\Psi} \varepsilon,1+2\tilde{L}_{\Psi} \varepsilon]\subset [1/2,3/2]$. Moreover, the matrix $\D{x}\Phi(x,\tau)$ is invertible and its inverse coincides with its pseudo-inverse. From the singular value decomposition of $\D{x}\Phi(x,\tau)$, we have that its pseudo-inverse $\D{x}\Phi(x,\tau)^\dagger$ is given by
        \begin{align*}
            \D{x}\Phi(x,\tau)^{\dagger} = U(x,\tau)\Sigma(x,\tau)^{\dagger} V(x,\tau)^\top,
        \end{align*}
        However, $\Sigma(x,\tau)^{\dagger}$ is simply the inverse of $\Sigma(x,\tau)$ which is well-defined since $\Sigma(x,\tau)$ is a diagonal matrix whose diagonal entries belong to the compact interval $[1-2\tilde{L}_{\Psi} \varepsilon,1+2\tilde{L}_{\Psi} \varepsilon]\subset [1/2,3/2]$. Therefore, we have that
        \begin{align*}
            \D{x}\Phi(x,\tau)^{-1}=\D{x}\Phi(x,\tau)^{\dagger} = U(x,\tau)\Sigma(x,\tau)^{-1} V(x,\tau)^\top,
        \end{align*}
        and, using the properties of the operator norm of matrices, we have that
        \begin{align*}
            \left|\D{x}\Phi(x,\tau)\right|&\leq |U(x,\tau)|\left|\Sigma(x,\tau)\right| |V(x,\tau)|\\
            \left|\D{x}\Phi(x,\tau)^{-1}\right|&\leq |U(x,\tau)|\left|\Sigma(x,\tau)^{-1}\right| |V(x,\tau)|
        \end{align*}
        Since $U(x,\tau)$ and $V(x,\tau)$ are orthonormal matrices, it follows that $|U(x,\tau)|=|V(x,\tau)|=1$. In addition, since $\Sigma(x,\tau)$ is a diagonal matrix whose diagonal entries belong to the compact interval $[1-2\tilde{L}_{\Psi} \varepsilon,1+2\tilde{L}_{\Psi} \varepsilon]\subset [1/2,3/2]$, we have that $|\Sigma(x,\tau)|\leq 1+2\tilde{L}_{\Psi} \varepsilon\leq \frac{3}{2}$, and $\left|\Sigma(x,\tau)^{-1}\right|\leq \frac{1}{1-2\tilde{L}_{\Psi} \varepsilon}\leq 2$. However, since $0\leq \varepsilon\leq \frac{1}{4\tilde{L}_{\Psi}}$, then $\frac{1}{1-2\tilde{L}_{\Psi} \varepsilon}\leq1+4\tilde{L}_{\Psi}\varepsilon$. Therefore, we have that
        \begin{align*}
            \left|\D{x}\Phi(x,\tau)\right|&\leq 1+2\tilde{L}_{\Psi} \varepsilon, &
            \left|\D{x}\Phi(x,\tau)^{-1}\right|&\leq1+4\tilde{L}_{\Psi}\varepsilon.
        \end{align*}
        It follows that the inverse of the Jacobian matrix $\D{}\Psi(x,\tau)$ is well-defined and is given by \eqref{inversetransfor}, \cite[p.146]{bernstein2009matrix}. The proof of the claim is concluded by defining $L_\Psi:=4 \tilde{L}_{\Psi}$. \hfill $\blacksquare$\\
\subsection{Proof of \thref{claim:boundedness_of_g}}
\textbf{Proof:} The map $g$ from item \ref{prop:near_identity_transformation_5} in \thref{prop:near_identity_transformation} has the explicit form
        \begin{align*}
            g(x,\tau,\varepsilon)&=  {\textstyle\sum}_{k=1}^5 G_i(x,\tau,\varepsilon)\,g_i(x,\tau,\varepsilon),
        \end{align*}
        where the matrix-valued maps $G_i$ are uniformly bounded and the maps $g_i$ are given by
        \begin{align*}
            g_1(x,\tau)&={v}_1\circ\Psi^{-1}(x,\tau), & g_2(x,\tau)&= -f_2\circ\Psi^{-1}(x,\tau), \\
            g_3(x,\tau)&=-f_1\circ\Psi^{-1}(x,\tau), & g_4(x,\tau)&={v}_2\circ\Psi^{-1}(x,\tau), \\
            g_5(x,\tau)&= -f_2\circ\Psi^{-1}(x,\tau).
        \end{align*}
        In this case, since $\delta_1=\delta_2=0$, we have that $\hat{f}_k(x,\tau)=f_k(x,\tau)$, and therefore we obtain that
        \begin{align*}
            f_1(x,\tau)&=\textstyle\sum_{i=1}^r\sum_{j=1}^2 b_{i,j}\sqrt{2\omega_i}\xi_i(J(x)+\omega_i\tau), & f_2(x,\tau)&=b_0(x),
        \end{align*}
        where $\xi_1(s)=\cos(s)$ and $\xi_2(s)=\sin(s)$. Clearly, for all $(x,\tau)\in\mathbb{R}^n\times\mathbb{R}_{\geq0}$, we have that
        \begin{align*}
            |f_1(x,\tau)|&\leq \textstyle\sum_{i=1}^r\sum_{j=1}^2\sqrt{2\omega_i}|b_{i,j}|, & |f_2(x,\tau)|&\leq |\nabla J(x)|\leq M_J.
        \end{align*}
        where used item b) in \thref{asmp:cost_function_conditions_es_1}. In addition, direct differentiation shows that
        \begin{align*}
            |\D{x} f_1(x,\tau)|&\leq \textstyle\sum_{i=1}^r\sum_{j=1}^2\sqrt{2\omega_i}|b_{i,j}||\nabla J(x)|,
        \end{align*}
        which is also uniformly bounded due to item b) in \thref{asmp:cost_function_conditions_es_1}. Since $v_1$ is the integral of $f_1$ with respect to $\tau$ and is periodic in $\tau$, it follows that $v_1$ is also uniformly bounded. Similarly, $v_2$ is the integral with respect to $\tau$ of terms that (smoothly) depend on $f_1$, $\D{x} f_1$, and $f_2$, all of which are uniformly bounded, and is periodic in $\tau$. It follows that $v_2$ is also uniformly bounded. Finally, since $\Psi^{-1}$ is $\mathcal{C}^1$-diffeomorphism, it follows that all the maps $g_i$ for $i\in\{1,\dots,5\}$ are uniformly bounded. Therefore, the remainder map $g$ is also uniformly bounded. \hfill $\blacksquare$
\end{document}